\documentclass[11pt]{amsart}
\usepackage{amscd,amssymb}
\usepackage{amsthm,amsmath,amssymb}
\usepackage[matrix,arrow]{xy}

\theoremstyle{definition}
\newtheorem{example}[equation]{Example}
\newtheorem{definition}[equation]{Definition}
\newtheorem{theorem}[equation]{Theorem}
\newtheorem{lemma}[equation]{Lemma}
\newtheorem{corollary}[equation]{Corollary}
\newtheorem{conjecture}[equation]{Conjecture}
\newtheorem{question}[equation]{Question}
\newtheorem*{question*}{Question}

\newtheorem*{problem*}{Problem}

\theoremstyle{remark}
\newtheorem{remark}[equation]{Remark}

\makeatletter\@addtoreset{equation}{section} \makeatother

\def\C {\mathbb{C}}
\def\Z {\mathbb{Z}}
\def\Q {\mathbb{Q}}
\def\P {\mathbb{P}}
\def\F {\mathbb{F}}
\def\GL {\mathrm{GL}}
\def\PGL {\mathrm{PGL}}
\def\SL {\mathrm{SL}}
\def\PSL {\mathrm{PSL}}
\def\A {\mathrm{A}}
\def\SS {\mathrm{S}}
\def\PSp {\mathrm{PSp}_4(\F_3)}
\def\Sp {\mathrm{Sp}_4(\F_3)}
\def\lct {\mathrm{lct}}
\def\qlin {\sim_{\Q}}

\def\le {\leqslant}
\def\ge {\geqslant}

\author{Ivan Cheltsov and Constantin Shramov}

\title{On exceptional quotient singularities}

\thanks{The authors were partially supported by AG Laboratory GU-HSE, RF
government grant 11~11.G34.31.0023. The first author was partially
supported by the~grants NSF DMS-0701465 and EPSRC EP/E048412/1.
The~second~author was partially supported by the~grants
RFFI~08-01-00395-a, RFFI~11-01-00185-a, RFFI~11-01-00336-a,
N.Sh.-1987.2008.1, N.Sh.-4713.2010.1 and EPSRC EP/E048412/1.}

\address{University of Edinburgh, Edinburgh EH9 3JZ, UK, \texttt{I.Cheltsov@ed.ac.uk}}

\address{Steklov Institute of Mathematics, Moscow 119991, Russia, \texttt{shramov@mccme.ru}}

\address{Laboratory of Algebraic Geometry, GU-HSE, 7 Vavilova street, Moscow 117312, Russia}%

\begin{document}

\begin{abstract}
We study  exceptional quotient singularities. In particular, we
prove an exceptionality criterion in terms of
the~$\alpha$-invariant of Tian, and utilize it to classify
four-dimensional and five-dimensional exceptional quotient
singularities.
\end{abstract}

\maketitle

\tableofcontents

We assume that all varieties are projective, normal, and defined
over $\mathbb{C}$.

\section{Introduction}
\label{section:intro}

Let $X$ be a~smooth Fano variety (see \cite{IsPr99}) of dimension $n$, let
$g=g_{i\overline{j}}$ be a~K\"ahler metric with a K\"ahler form
$$
\omega=\frac{\sqrt{-1}}{2\pi}\sum g_{i\overline{j}}dz_i\wedge d\overline{z}_j\in \mathrm{c}_1\big(X\big).%
$$

\begin{definition}
\label{definition:KE} The metric $g$ is a~K\"ahler--Einstein
metric if $\mathrm{Ric}\big(\omega\big)=\omega$, where
$\mathrm{Ric}(\omega)$ is a~Ricci curvature of the~metric $g$.
\end{definition}

Let $\bar{G}\subset\mathrm{Aut}(X)$ be a~compact subgroup. Suppose
that $g$ is $\bar{G}$-invariant.

\begin{definition}
\label{definition:alpha-invariant} Let $P_{\bar{G}}(X,g)$ be the~
set of $C^2$-smooth $\bar{G}$-invariant functions $\varphi$ such
that
$$
\omega+\frac{\sqrt{-1}}{2\pi}\partial\overline{\partial}\varphi>0
$$
and $\sup_X \varphi =0$. Then the~$\bar{G}$-invariant
$\alpha$-invariant of the~variety $X$ is the~number
$$
\alpha_{\bar{G}}(X)=\mathrm{sup}\Bigg\{\lambda\in\mathbb{Q}\ \Big\vert\ \exists\ C\in\mathbb{R}\ \text{such that}\ \int_X e^{-\lambda\varphi}\omega^n\leqslant C\ \text{for any}\ \varphi\in P_{\bar{G}}\big(X,g\big)\Bigg\}.%
$$
\end{definition}

The number $\alpha_{\bar{G}}(X)$ was introduced in \cite{Ti87} and
\cite{TiYa87} and now it is  called the $\alpha$-invariant of
Tian.

\begin{theorem}[{\cite{Ti87}}]
\label{theorem:KE} The Fano variety $X$ admits
a~$\bar{G}$-invariant K\"ahler--Einstein metric~if
$\alpha_{\bar{G}}(X)>n/(n+1)$.
\end{theorem}

The normalized K\"ahler--Ricci flow on the~smooth Fano $X$ is
defined by the~equation
\begin{equation}
\label{equation:flow} \left\{\aligned
&\frac{\partial \omega(t)}{\partial t}=-\mathrm{Ric}\big(\omega(t)\big)+\omega(t),\\
&\omega(0)=\omega,\\
\endaligned
\right.
\end{equation}
where $\omega(t)$ is a~K\"ahler form such that
$\omega(t)\in\mathrm{c}_1(X)$, and $t\in\mathbb{R}_{{}\geqslant
0}$. It follows from~\cite{Cao85} that the~solution $\omega(t)$ to
(\ref{equation:flow}) exists for every $t>0$.

\begin{theorem}[\cite{TiZhu07}]
\label{theorem:Perelman} If $X$ admits a~K\"ahler--Einstein metric
with a~K\"ahler form~$\omega_{KE}$, then any~solution to
(\ref{equation:flow}) converges to $\omega_{KE}$ in the~sense~of~Cheeger--Gromov.%
\end{theorem}

The normalized K\"ahler--Ricci iteration on the~smooth Fano variety~$X$ is
defined by the~equation
\begin{equation}
\label{equation:iteration} \left\{\aligned
&\omega_{i-1}=\mathrm{Ric}\big(\omega_{i}\big),\\
&\omega_{0}=\omega,\\
\endaligned
\right.
\end{equation}
where $\omega_{i}$ is a~K\"ahler form such that
$\omega_{i}\in\mathrm{c}_1(X)$. It follows from \cite{Yau78} that
the~solution~$\omega_{i}$
to (\ref{equation:iteration}) exists for every $i\geqslant 1$.%

\begin{theorem}[{\cite{Rub08}}]
\label{theorem:Yanir} If $\alpha_{\bar{G}}(X)>1$ then $X$ admits
a~$\bar{G}$-invariant K\"ahler--Einstein metric with a~K\"ahler
form $\omega_{KE}$ and any~solution to (\ref{equation:iteration})
converges to $\omega_{KE}$ in $C^{\infty}(X)$-topology.
\end{theorem}

Smooth Fano varieties that satisfy all hypotheses of
Theorem~\ref{theorem:Yanir} do exist.

\begin{example}
\label{example:dP5} Let $X$ be a~smooth del Pezzo surface such
that $K_{X}^{2}=5$. Then $X$ is unique and
$\mathrm{Aut}(X)\cong\SS_{5}$. Moreover, one can show that
$\alpha_{\bar{G}}(X)=2$  in the case when
$\bar{G}\cong\SS_{5}$ or $\bar{G}\cong\A_{5}$ (see
\cite[Example~1.11]{Ch07b} and \cite[Theorem~A.3]{ChSh08c}).
\end{example}

Suppose now that $X=\mathbb{P}^{n}$ (the simplest possible case).
Then the~Fubini--Studi metric on $\mathbb{P}^{n}$ is
K\"ahler--Einstein. Moveover, if $\bar{G}$ is the~maximal compact
subgroup of $\mathrm{Aut}(\mathbb{P}^{n})$, then the only
$\bar{G}$-invariant metric on $\mathbb{P}^{n}$ is
the~Fubini--Studi metric and we have
$\alpha_{\bar{G}}(\mathbb{P}^{n})=+\infty$ by
Definition~\ref{definition:alpha-invariant}. In particular,
the~solution to (\ref{equation:iteration}) is trivial (and
constant) in the latter case, since the~initial metric $g$ must be
the~Fubini--Studi metric. On the other hand,  the convergence of
any~solution to (\ref{equation:iteration}) is not clear in the
case when $\bar{G}$ is a finite group. So, Yanir Rubinstein asked
the~following question in the~Spring 2009.

\begin{question}
\label{question:Yanir} Is there a~finite subgroup
$\bar{G}\subset\mathrm{Aut}(\mathbb{P}^{n})$ such that
$\alpha_{\bar{G}}(\mathbb{P}^{n})>1$?
\end{question}

This paper is inspired by Question~\ref{question:Yanir}. In
particular, we will show that the answer to
Question~\ref{question:Yanir} is positive in the case
when~$n\leqslant 4$, which follows from
\cite[Theorem~A.3]{ChSh08c} and
Theorems~\ref{theorem:Shokurov-n-2}, \ref{theorem:Dima-Yura-n-3},
\ref{theorem:Vanya-Kostya-n-4}, \ref{theorem:Vanya-Kostya-n-5} and
\ref{theorem:Dima-Yura-strong}.

\medskip

It came as a surprise that Question~\ref{question:Yanir} is
strongly related to the notion of exceptional singularity that was
introduced by Vyacheslav Shokurov in \cite{Sho00}. Let us recall
this notion. Let $(V\ni O)$ be a~germ of Kawamata log terminal
singularity (see \cite[Definition~3.5]{Ko97}).

\begin{definition}[{\cite[Definition~1.5]{Sho00}}]
\label{definition:exceptional} The singularity $(V\ni O)$ is~said
to be \emph{exceptional} if for every effective
$\mathbb{Q}$-divisor $D_{V}$ on the~variety $V$ such that
$(V,D_{V})$ is log canonical (see \cite[Definition~3.5]{Ko97}) and
for every resolution of singularities $\pi\colon U\to V$ there
exists at most one $\pi$-exceptional divisor $E\subset
U$~such~that $a(V,D_{V},E)=-1$, where the rational number
$a(V,D_{V},E)$ can be defined through the equivalence
$$
K_{U}+D_{U}\sim_{\mathbb{Q}}\pi^{*}\Big(K_{V}+D_{V}\Big)+
\sum a\big(V,D_{V},E\big)E,%
$$
where the sum is taken over all $f$-exceptional divisors, and
$D_{U}$ is the~proper transform of the~divisor $D_{V}$ on
the~variety $U$.
\end{definition}

One can show that exceptional Kawamata log terminal
singularities are straightforward generalizations of the Du Val
singularities of type $\mathbb{E}_{6}$, $\mathbb{E}_{7}$ and
$\mathbb{E}_{8}$ (cf. Theorem~\ref{theorem:Shokurov-n-2}), which
partially justifies the word ``exceptional'' in
Definition~\ref{definition:exceptional}.

\begin{remark}
\label{remark:smooth-points} One can easily check (for example, by
applying Theorem~\ref{theorem:exceptional-criterion}) that the
singularity $(V\ni O)$ is not exceptional if $V$ is smooth and
$\mathrm{dim}(V)\geqslant 2$.
\end{remark}

It follows from~\cite{Sho93},~\cite{IshiiPr01} and \cite{MarPr99}
that exceptional Kawamata log terminal singularities do~exist
in~dimensions~$2$~and~$3$. The existence in dimension $4$ follows
from \cite{JoKo01} and~\cite[Theorem~4.9]{Pr98plt}. Actually,
exceptional Kawamata log terminal singularities exist in every
dimension (see Example~\ref{example:Kollar-Boyer}). We will see
later (cf. Theorem~\ref{theorem:alpha},
Remark~\ref{remark:reflections}, Theorem~\ref{theorem:criterion}
and Conjecture~\ref{conjecture:Vanya-Kostya-I}) that
Question~\ref{question:Yanir} is \emph{almost} equivalent to the
following

\begin{question}
\label{question:exceptional} Are there exceptional \emph{quotient}
singularities of dimension $n+1$?
\end{question}

Recall that quotient singularities are always Kawamata
log~terminal by \cite[Proposition~3.16]{Ko97}. So
Question~\ref{question:exceptional} fits well to
Definition~\ref{definition:exceptional}. Moreover, it follows
from~\cite{Sho00} and~\cite{MarPr99} that the answer to
Question~\ref{question:exceptional} is positive for $n=1$ and
$n=2$, respectively. The~purpose of this paper is to study
exceptional \emph{quotient} singularities and, in particular, to
give positive answers to Questions~\ref{question:Yanir} and
\ref{question:exceptional} for every~\mbox{$n\leqslant 4$}. In a
subsequent paper  we will show that the answers to
Questions~\ref{question:Yanir}
and~\ref{question:exceptional} are still positive for $n=5$ and are
surprisingly negative for $n=6$ (see \cite{ChSh10}).
So it is hard to predict what would be the answer
to Question~\ref{question:Yanir} in general.
However, we
still believe in the following

\begin{conjecture}
\label{conjecture:Vanya-Kostya} For every $N\in\Z_{{}>0}$ there exist
exceptional quotient singularities of dimension greater than~$N$.
\end{conjecture}

Let $G$ be a~finite subgroup in $\GL_{n+1}(\mathbb{C})$,
where $n\geqslant 1$. Denote by~$Z(G)$ and~\mbox{$[G,G]$} the~center and
the~commutator of group $G$, respectively. Let
\mbox{$\phi\colon\GL_{n+1}(\mathbb{C})\to\mathrm{Aut}(\mathbb{P}^{n})\cong
\PGL_{n+1}(\mathbb{C})$}
be the~natural projection. Put  $\bar{G}=\phi(G)$ and put
$$
\mathrm{lct}\Big(\mathbb{P}^{n},\bar{G}\Big)=\mathrm{sup}\left\{\lambda\in\mathbb{Q}\ \left|%
\aligned
&\text{the~log pair}\ \left(\mathbb{P}^{n}, \lambda D\right)\ \text{has log canonical singularities}\\
&\text{for every $\bar{G}$-invariant effective $\mathbb{Q}$-divisor}\ D\sim_{\mathbb{Q}} -K_{\mathbb{P}^{n}}\\
\endaligned\right.\right\}.%
$$

\begin{theorem}[{see e.\,g. \cite[Theorem~A.3]{ChSh08c}}]
\label{theorem:alpha} One has
$\mathrm{lct}(\mathbb{P}^{n},\bar{G})=\alpha_{\bar{G}}(\mathbb{P}^{n})$.
\end{theorem}

The number $\mathrm{lct}(\mathbb{P}^{n},\bar{G})$ is usually
called $\bar{G}$-equivariant \emph{global log canonical threshold}
of $\mathbb{P}^{n}$. Despite the fact that
$\mathrm{lct}(\mathbb{P}^{n},\bar{G})=\alpha_{\bar{G}}(\mathbb{P}^{n})$,
we still prefer to work with the number
$\mathrm{lct}(\mathbb{P}^{n},\bar{G})$ throughout this paper,
because it is easier to
handle than~$\alpha_{\bar{G}}(\mathbb{P}^{n})$. For
example, it follows immediately from
Definition~\ref{definition:G-threshold} that
$\mathrm{lct}(\mathbb{P}^{n},\bar{G})\leqslant d/(n+1)$ if
the~group $G$ has a~semi-invariant of degree~$d$ (a semi-invariant
of the group $G$ is a polynomial whose zeroes define a
$\bar{G}$-invariant hypersurface in $\mathbb{P}^{n}$).

\begin{remark}
\label{remark:semiinvariants-vs-invariants} A semi-invariant of
the~group $G$ is its invariant if $Z(G)\subseteq[G,G]$
and~$\bar{G}$ is a non-abelian simple group.
\end{remark}

Recall that an element $g\in G$ is
called a \emph{reflection} (or sometimes a \emph{quasi-reflection})
if there is
a hyperplane in $\mathbb{P}^{n}$ that is pointwise fixed by
$\phi(g)$ (cf. \cite[\S 4.1]{Spr77}).

\begin{remark}
\label{remark:reflections} Let $R\subseteq G$ be a~subgroup
generated by all reflections. Then
the~quotient $\mathbb{C}^{n+1}\slash R$ is isomorphic
to $\mathbb{C}^{n+1}$ (see \cite{SheTo54},
\cite[Theorem~4.2.5]{Spr77}).
Moreover, the~subgroup $R\subseteq G$ is
normal, and the singularity $\C^{n+1}/G$ is isomorphic to the
singularity $\C^{n+1}/(G/R)$.
Note that the~subgroup $R$ is
trivial if $G\subset\SL_{n+1}(\mathbb{C})$.
If~$G$ is a~trivial group, then the singularity $\C^{n+1}/G\cong\C^{n+1}$
is not exceptional by
Remark~\ref{remark:smooth-points}.
\end{remark}

Thus to answer Question~\ref{question:exceptional} one can always
assume that the group $G$ does not contain reflections. On the
other hand, one can easily check that there exists a~finite
subgroup $G^{\prime}\subset\SL_{n+1}(\mathbb{C})$ such~that
$\phi(G^{\prime})=\bar{G}$. So to answer
Question~\ref{question:Yanir} one can also assume that
$G\subset\SL_{n+1}(\mathbb{C})$, which implies, in
particular, that the group~$G$ does not contain reflections.
Moreover, if the group $G$ does not contain reflections, then the
singularity $\mathbb{C}^{n+1}\slash G$ is exceptional if and only
if the~singularity $\mathbb{C}^{n+1}\slash G^{\prime}$ is
exceptional thanks to the following

\begin{theorem}
\label{theorem:criterion} Let $G$ be a~finite subgroup in
$\GL_{n+1}(\mathbb{C})$ that does not contain reflections.
Then
\begin{itemize}
\item the~singularity $\C^{n+1}/G$ is exceptional if $\mathrm{lct}(\mathbb{P}^{n},\bar{G})>1$,%
\item the~singularity $\C^{n+1}/G$ is not exceptional if either $\mathrm{lct}(\mathbb{P}^{n},\bar{G})<1$ or $G$ has a~semi-invariant of degree~at~most~$n+1$,%
\item for any~subgroup
$G^{\prime}\subset\GL_{n+1}(\mathbb{C})$ such that $G^{\prime}$
does not contain reflections and
$\phi(G^{\prime})=\bar{G}$, the~singularity $\C^{n+1}/G$ is
exceptional if and only if the~singularity $\mathbb{C}^{n+1}\slash G^{\prime}$ is exceptional.%
\end{itemize}
\end{theorem}

\begin{proof}
All required assertions immediately follow from
Theorem~\ref{theorem:exceptional-quotient} (cf.
\cite[Proposition~3.1]{Pr00}, \cite[Lemma~3.1]{Pr00}).
\end{proof}

It should be pointed out that the assumption that $G$ contains no
reflections is crucial for Theorem~\ref{theorem:criterion}.

\begin{example}
\label{example:reflections} Let $G$ be a finite subgroup in
$\GL_4(\mathbb{C})$ that is the~subgroup number~$32$ in
\cite[Table~VII]{SheTo54}. Then the~group~$G$ is generated by
reflections (see \cite{SheTo54}), so that the~singularity
$\mathbb{C}^{4}\slash G$ is not exceptional by
Remark~\ref{remark:reflections}. On the other hand, it follows
from Theorem~\ref{theorem:Vanya-Kostya-n-4} that
$\mathrm{lct}(\mathbb{P}^{3},\bar{G})\geqslant 5/4$, because
$\bar{G}\cong\PSp$. It follows from
Theorem~\ref{theorem:Vanya-Kostya-n-4} that there exists a
subgroup $G^{\prime}\subset\SL_4(\mathbb{C})$ such that
$\bar{G}=\phi(G^{\prime})$ and the singularity
$\mathbb{C}^{4}\slash G^{\prime}$ is exceptional.
One can produce similar examples for two-dimensional and three-dimensional
singularities.
\end{example}

By Theorem~\ref{theorem:criterion} and \cite[\S\,4.5]{Spr77}, if
$G$ is a finite subgroup in $\GL_2(\C)$ that does not contain
reflections, then the~singularity $\C^2/G$ is exceptional if and
only if $G$ has no semi-invariants of degree at most $2$. A
similar result holds in dimension $3$.

\begin{theorem}[{\cite[Theorem~1.2]{MarPr99}}]
\label{theorem:Dima-Yura} Let $G$ be a finite group in $\GL_3(\C)$
that does not contain reflections. Then the~singularity $\C^3/G$
is exceptional if and only if $G$ does not have semi-invariants of
degree at most $3$.
\end{theorem}

For finite subgroups in $\GL_4(\C)$, the~assertion of
Theorem~\ref{theorem:Dima-Yura} is no longer true.

\begin{example}[{\cite[Example~3.1]{Pr00}}]
\label{example:dim-4-imprimitive} Let
$\Gamma\subset\SL_2(\mathbb{C})$ be a~binary icosahedron
group. Put
$$
G=\left\{\aligned\left.\aligned\left(
\begin{array}{cccc}
a_{11} & a_{12} & 0&0 \\
a_{21} & a_{22} & 0&0 \\
0& 0& b_{11} & b_{12} \\
0& 0& b_{21} & b_{22} \\
\end{array}
\right)\ \endaligned\right|\ \left(\begin{array}{cc}
a_{11} & a_{12} \\
a_{21} & a_{22} \\
\end{array}
\right)\in \Gamma\ni
\left(\begin{array}{cc}
b_{11} & b_{12} \\
b_{21} & b_{22} \\
\end{array}
\right)\endaligned\right\}\subset\SL_4\Big(\mathbb{C}\Big),
$$
where $a_{ij}\in\mathbb{C}\ni b_{ij}$. Then $G$ does not have
semi-invariants of degree at most $4$, because $\Gamma$ does not
have semi-invariants of degree at most $4$ (see \cite[\S
4.5]{Spr77}). On the other hand, it follows from
\cite[Proposition~2.1]{Pr00} that the~singularity $\C^4/G$~is~not
exceptional (cf. Corollary~\ref{corollary:primitive}).
\end{example}

Actually, it is possible to modify the assertion of
Theorem~\ref{theorem:Dima-Yura} so that its new version can be
generalized to higher dimensions.

\begin{definition}[{\cite{Bli17}}]
\label{definition:primitive} The subgroup
$G\subset\GL_{n+1}(\mathbb{C})$ is said to be primitive if
there is~no~non-trivial decomposition
$\mathbb{C}^{n+1}=\bigoplus_{i=1}^{r}V_{i}$ such that
for any $g\in G$ and any $i$ there is some $j=j(g)$ such that
$g(V_{i})=V_{j}$.
\end{definition}

If $G$ is primitive, then $\bar{G}\cong G\slash Z(G)$ by Schur's
lemma. It follows from \cite[Proposition~2.1]{Pr00} that $G$ must
be primitive if $\C^{n+1}/G$ is exceptional (we give a short proof
of this fact in Corollary~\ref{corollary:primitive}). Moreover,
primitivity plays a crucial role in the main result of this
paper, which is the~following

\begin{theorem}
\label{theorem:Vanya-Kostya-invariants}  Let $G$ be a finite
subgroup in $\GL_{n+1}(\C)$ that does not contain reflections.
Suppose that $n\le 4$. Then the~following conditions are
equivalent:
\begin{itemize}
\item the~singularity $\C^{n+1}/G$ is exceptional,%
\item $\lct(\P^n, \bar{G})\ge (n+2)/(n+1)$,
\item the~group $G$ is primitive and has no semi-invariants
of degree at most $n+1$.%
\end{itemize}
\end{theorem}

\begin{proof}
The required assertion follows from
Theorems~\ref{theorem:Dima-Yura},
\ref{theorem:exceptional-quotient}, \ref{theorem:primitive},
\ref{theorem:Dima-Yura-strong}, \ref{theorem:Vanya-Kostya-n-4}
and~\ref{theorem:Vanya-Kostya-n-5}.
\end{proof}

It appears that in higher dimensions exceptionality cannot be
expressed in terms of primitivity and absence of semi-invariants
of small degree. In particular, there are non-exceptional
six-dimensional quotient singularities arising from primitive
subgroups without reflections in $\GL_6(\C)$ that have no
semi-invariants of degree at most~$6$ (see
Example~\ref{example:Severi}). On the other hand, it follows from
Theorem~\ref{theorem:Vanya-Kostya-invariants} that we may expect
the sufficient condition for exceptionality in
Theorem~\ref{theorem:criterion} to be a necessary condition as
well. Namely, inspired by
Theorem~\ref{theorem:Vanya-Kostya-invariants} and
\cite[Question~1]{Ti90b} we believe in the following

\begin{conjecture}
\label{conjecture:Vanya-Kostya-I} Let $G$ be a finite subgroup in
$\GL_{n+1}(\C)$ that does not contain reflections. Then
the~singularity $\mathbb{C}^{n+1}/G$ is exceptional if and only if
$\mathrm{lct}(\mathbb{P}^{n},\bar{G})>1$.
\end{conjecture}

It follows from Theorem~\ref{theorem:Vanya-Kostya-invariants} that
Conjecture~\ref{conjecture:Vanya-Kostya-I} holds for $n\leqslant
4$. In a subsequent paper  we will show that
Conjecture~\ref{conjecture:Vanya-Kostya-I} holds for $n=5$ and
$n=6$ (see \cite{ChSh10}). Note that
Conjecture~\ref{conjecture:Vanya-Kostya-I} is a special case of
Conjecture~\ref{conjecture:stabilization}.

\medskip

To apply Theorem~\ref{theorem:Vanya-Kostya-invariants} we may
assume that $G\subset\SL_{n+1}(\mathbb{C})$, since there
exists a~finite subgroup $G^{\prime}\subset\SL_{n+1}(\mathbb{C})$
such~that $\phi(G^{\prime})=\bar{G}$. On the other
hand, it is well known that there are at most finitely many
primitive finite subgroups in $\SL_{n+1}(\mathbb{C})$ up to
conjugation (see \cite{Col07}). Primitive finite subgroups of
$\SL_2(\mathbb{C})$ are group-theoretic counterparts of Platonic
solids and each of them gives rise to an exceptional singularity
(see Theorem~\ref{theorem:Shokurov-n-2}). Primitive finite
subgroups of $\SL_3(\mathbb{C})$ are classified by
Blichfeldt in \cite{Bli17}. Prokhorov and Markushevich used
Blichfeldt's classification in~\cite{MarPr99} to obtain an
explicit classification of the subgroups in
$\SL_3(\mathbb{C})$ corresponding to three-dimensional
exceptional quotient singularities (see
Theorem~\ref{theorem:Dima-Yura-n-3}). For dimension $2$ the same
was done by Shokurov (see Theorem~\ref{theorem:Shokurov-n-2}).
Similar classification is possible in dimensions~$4$ and~$5$, since
primitive finite subgroups of $\SL_4(\mathbb{C})$
and~\mbox{$\SL_5(\mathbb{C})$} are classified in \cite{Bli17} and
\cite{Br67}, respectively. In fact, we obtain a complete list of
finite subgroups in $\SL_4(\mathbb{C})$ and
$\SL_5(\mathbb{C})$ that satisfy all hypotheses of
Theorem~\ref{theorem:Vanya-Kostya-invariants} (see
Theorems~\ref{theorem:Vanya-Kostya-n-4}
and~\ref{theorem:Vanya-Kostya-n-5}).

\medskip

While the exceptionality of a quotient singularity
$\mathbb{C}^{n+1}/G$ depends on a lower bound for a global log
canonical threshold $\lct(\P^n, \bar{G})$, it is interesting to
find upper bounds for $\lct(\P^n, \bar{G})$ as well. Using
\cite[\S\,4.5]{Spr77}, \cite{YauYu93} and a bit of direct
computation, we see that it follows from
Corollary~\ref{corollary:primitive-baby} that
$$
\mathrm{lct}\big(\mathbb{P}^{n},\bar{G}\big)\leqslant\left\{\aligned
&6\ \text{if}\ n=1,\\
&2\ \text {if}\ n=2,\\
&3\ \text {if}\ n=3.\\
\endaligned
\right.
$$

\begin{theorem}
\label{theorem:thomas} The inequality
$\mathrm{lct}(\mathbb{P}^{n}, {\bar{G}})\leqslant 4(n+1)$ holds
for every $n\geqslant 1$. Moreover, if
$n\geqslant 23$, then
$\mathrm{lct}(\mathbb{P}^{n}, {\bar{G}})\leqslant 12(n+1)/5$.
\end{theorem}

\begin{proof}
Let $p$ be any~prime number which does not divide $|G|$. Then $G$
has a~semi-invariant of degree at most $(p-1)(n+1)$ by
\cite[Lemma~2]{Th81}. Thus, it follows from
Definition~\ref{definition:G-threshold} that
$\mathrm{lct}(\mathbb{P}^{n},\bar{G})\leqslant p-1$. On the other
hand, it follows from the~Bertrand's postulate (see \cite{Ram19})
that there is a~prime number $p^{\prime}$ such that
\mbox{$2n+3<p^{\prime}<2(2n+3)$}, which implies that $p^{\prime}\leqslant
4n+5$. If $G$ is not primitive, then
$\mathrm{lct}(\mathbb{P}^{n},\bar{G})\leqslant 1$ by
Corollary~\ref{corollary:primitive-baby}. If $G$ is is primitive,
then $p^{\prime}$ does not divide~\mbox{$|G|$}
by~\cite[Theorem~1]{FeTh63}, which completes the proof of the first
assertion of the theorem. A similar argument with an additional
use of~\cite{Na52} gives the second assertion for $n\ge 23$.
\end{proof}

In fact, we expect the~following to be true (cf. \cite{Th81}).

\begin{conjecture}
\label{conjecture:Thomas} There exists a universal constant
$C\in\mathbb{R}$ such that $\mathrm{lct}(\mathbb{P}^{n},
{\bar{G}})\leqslant C$ for any finite subgroup
$\bar{G}\subset\mathrm{Aut}(\mathbb{P}^{n})$~and~for any
$n\geqslant 1$.
\end{conjecture}

\bigskip

Let us describe the structure of the paper. In
Section~\ref{section:preliminaries} we collect auxiliary results.
In Section~\ref{section:plt-blow-up} we prove the~exceptionality
criterion for a singularity $\C^{n+1}/G$. In
Sections~\ref{section:4} we classify exceptional quotient
singularities in dimension~$4$ (see
Theorem~\ref{theorem:Vanya-Kostya-n-4}). In
Sections~\ref{section:5} we classify exceptional quotient
singularities in dimension~$5$ (see
Theorem~\ref{theorem:Vanya-Kostya-n-5}). In
Appendix~\ref{section:HM} we prove
Corollary~\ref{corollary:Heisenberg-representations} and
Theorem~\ref{theorem:HM-splitting} that are used in
Section~\ref{section:5}.

Many of our results can be obtained by direct computations using
\cite{Atlas}.

Throughout the~paper we use the~following standard notation:
the~symbol $\mathbb{Z}_n$ denotes the~cyclic group of order $n$,
the~symbol $\mathbb{F}_n$ denotes the~finite field consisting of
$n$ elements, the~symbol $\SS_n$ denotes the~symmetric
group of degree $n$, the~symbol~$\A_n$ denotes
the~alternating group of degree $n$, the~symbols $\GL$,
$\PGL$, $\SL$, $\PSL$, $\Sp$ and
$\PSp$ denote the~corresponding algebraic groups. The~symbol~$k.G$
denotes a~central extension of a~group $G$ with the~center
$\mathbb{Z}_k$ (this might be non-unique).

We would like to thank I.\,Arzhantsev, S.\,Galkin, V.\,Dotsenko,
A.\,Khoroshkin,
S.\,Lok\-tev,~D.\,Pasech\-nik,~V.\,Po\-pov,~Yu.\,Pro\-kho\-rov,
S.\,Rybakov, L.\,Rybnikov and V.\,Vologodsky for very useful and
fruitful discussions. We thank T.\,K\"oppe for helping us access
the~classical German literature on Invariant Theory.

\section{Preliminaries}%
\label{section:preliminaries}

Throughout the paper we use the standard language of the
singularities of pairs (see \cite{Ko97}). By strictly log
canonical singularities we mean log canonical singularities that
are not Kawamata log terminal (see \cite[Definition~3.5]{Ko97}).

Let $X$ be a~variety, let $B_{X}$ and $D_{X}$ be effective
$\mathbb{Q}$-divisors on the~variety $X$~such that
the~singularities of the~log pair $(X,B_{X})$~are Kawamata log
terminal, and \mbox{$K_{X}+B_{X}+D_{X}$} is a~$\mathbb{Q}$-Cartier
divisor. Let $Z\subseteq X$ be a~closed non-empty  subvariety.

\begin{definition}
\label{definition:lct-local} The log canonical threshold of
the~boundary $D_{X}$ along $Z$ is
$$
\mathrm{c}_{Z}\big(X,B_{X},D_{X}\big)=\mathrm{sup}\left\{\lambda\in\mathbb{Q}\
\Big|\ \text{the~pair}\
 \big(X, B_{X}+\lambda D_{X}\big)\ \text{is log canonical along}~Z\right\}.%
$$
\end{definition}

Note that the~log pair $(X,B_{X}+D_{X})$ is Kawamata log terminal
along $Z$ if and only if $\mathrm{c}_{Z}(X,B_{X},D_{X})>1$. For
simplicity, we put
$\mathrm{c}(X,B_{X},D_{X})=\mathrm{c}_{X}(X,B_{X},D_{X})$. We put
$\mathrm{c}_{Z}(X,D_{X})=\mathrm{c}_Z(X,B_{X},D_{X})$ in the~case
when $B_{X}=0$. For simplicity, we also put
$\mathrm{c}(X,D_{X})=\mathrm{c}_X(X,D_{X})$.

\medskip

Apart from some rare but important occasions (cf.
Section~\ref{section:plt-blow-up}), we only need to consider the
case when $B_{X}=0$. So from now on we assume that $B_{X}=0$.

Let $\pi\colon\bar{X}\to X$ be a~birational morphism such that $\bar{X}$ is
smooth. Then
$$
K_{\bar{X}}+D_{\bar{X}}\sim_{\mathbb{Q}}\pi^{*}\Big(K_{X}+D_{X}\Big)+\sum_{i=1}^{m}d_{i}E_{i},
$$
where $D_{\bar{X}}$ is a proper transform of the~divisor $D_{X}$
on the~variety $\bar{X}$, $d_{i}\in\mathbb{Q}$, and $E_{i}$ is an
exceptional divisor of the~morphism $\pi$. Put
$D_{\bar{X}}=\sum_{i=1}^{r}a_{i}\bar{D}_{i}$, where
$a_{i}\in\mathbb{Q}_{{}\geqslant 0}$, and $\bar{D}_{i}$ is a prime
Weil divisor on $\bar{X}$. Suppose that
$\sum_{i=1}^{r}\bar{D}_{i}+\sum_{i=1}^{m}E_{i}$ is a~divisor with
simple normal crossing. Put
$$
\mathcal{I}\Big(X, D_{X}\Big)=\pi_{*}\mathcal{O}_{\bar{X}}
\Bigg(\sum_{i=1}^{m}\lceil d_{i}\rceil E_{i}-\sum_{i=1}^{r}\lfloor a_{i}\rfloor \bar{D}_{i}\Bigg),%
$$
and let $\mathcal{L}(X, D_{X})$ be a~subscheme that corresponds to
the ideal sheaf \mbox{$\mathcal{I}(X, D_{X})$} (the sheaf
$\mathcal{I}(X, D_{X})$ is an ideal sheaf, because $D_{X}$ is an
effective divisor). Put $\mathrm{LCS}(X,
D_{X})=\mathrm{Supp}(\mathcal{L}(X, D_{X}))$.

\begin{remark}
\label{remark:log-canonical-subscheme} If $(X, D_{X})$ is log canonical, then $\mathcal{L}(X,D_{X})$ is reduced. %
\end{remark}

The subscheme $\mathcal{L}(X, D_{X})$ and the locus
$\mathrm{LCS}(X, D_{X})$ were introduced in~\cite{Sho93}. They are
called are called the~subscheme of log canonical singularities of
the~log pair $(X, D_{X})$ and the~locus of log~canonical
singularities of the~log pair $(X, D_{X})$, respectively. Note
that the ideal sheaf $\mathcal{I}(X, D_{X})$ is also known as
the~multiplier ideal sheaf of the~log pair $(X, D_{X})$ (see
\cite{La04}).

\begin{theorem}[{\cite[Theorem~9.4.8]{La04}}]
\label{theorem:Shokurov-vanishing} Let $H$ be a~nef and big
$\mathbb{Q}$-divisor on $X$ such that $K_{X}+D_{X}+H\equiv D$ for
some Cartier divisor $D$ on the~variety $X$. Then
$$
H^{i}\Bigg(\mathcal{I}\Big(X, D_{X}\Big)\otimes D\Bigg)=0
$$
for every $i\geqslant 1$.
\end{theorem}

\begin{corollary}[{\cite[Lemma~5.7]{Sho93}}]
\label{corollary:connectedness} Suppose that $-(K_{X}+D_{X})$ is
nef and big. Then the locus $\mathrm{LCS}(X, D_{X})$ is connected.
\end{corollary}

Let $\mathbb{LCS}(X, D_{X})$ be the~set that consists of all
possible centers of log canonical singularities of the~log pair
$(X, D_{X})$ (see \cite[Definition~2.2]{ChSh08c}).

\begin{remark}
\label{remark:hyperplane-reduction} Let $\mathcal{H}$ be a~linear
system on the~variety $X$ that has  no base points. Put $Z\cap
H=\sum_{i=1}^{k}Z_{i}$, where $H$ is~a~general divisor in
$\mathcal{H}$, and $Z_{i}$ is an irreducible subvariety in $H$.
Then $Z\in\mathbb{LCS}(X, D_{X})$ if and only if all subvarieties
$Z_{1},\ldots,Z_{k}$ are contained in the set $\mathbb{LCS}(H,
D_{X}\vert_{H})$.
\end{remark}

If $Z\in\mathbb{LCS}(X, D_{X})$ and no proper subvariety of $Z$ is
contained in $\mathbb{LCS}(X, D_{X})$, then $Z$ is said to be a
\emph{minimal} center in $\mathbb{LCS}(X, D_{X})$ or
\emph{minimal} center of log canonical singularities of the~log
pair $(X, D_{X})$.

\begin{lemma}[{\cite[Proposition~1.5]{Kaw97}}]
\label{lemma:centers} Suppose that $Z\in\mathbb{LCS}(X, D_{X})$
and $(X, D_{X})$ is log canonical. Let $Z^{\prime}$ be a~center in
$\mathbb{LCS}(X, D_{X})$ such that $\varnothing\ne Z\cap
Z^{\prime}=\sum_{i=1}^{k}Z_{i}$, where $Z_{i}\subsetneq Z$ is an
irreducible subvariety. Then $Z_{i}\in\mathbb{LCS}(X, D_{X})$ for
every $i\in\{1,\ldots,k\}$.
\end{lemma}

\begin{theorem}[{\cite[Theorem~1]{Kaw98}}]
\label{theorem:Kawamata} Suppose that $Z\subset X$ is a~minimal
center in $\mathbb{LCS}(X, D_{X})$  and $(X, D_{X})$ is log
canonical. Then $Z$ is normal and has at most rational
singularities. Let $\Delta$ be an~ample $\mathbb{Q}$-Cartier
$\mathbb{Q}$-divisor on $X$. Then there exists an~effective
$\mathbb{Q}$-divisor $B_{Z}$ on the~variety $Z$ such that
$$
\Big(K_{X}+D_{X}+\Delta\Big)\Big\vert_{Z}\sim_{\mathbb{Q}} K_{Z}+B_{Z},%
$$
and $(Z,B_{Z})$ has Kawamata log terminal singularities.
\end{theorem}

Let $\bar{G}\subseteq\mathrm{Aut}(X)$ be a~finite subgroup such
that $D_{X}$ is $\bar{G}$-invariant. Then
$g(Z)\in\mathbb{LCS}(X,D_{X})$ for every $g\in\bar{G}$, and
the~locus $\mathrm{LCS}(X,D_{X})$ is $\bar{G}$-invariant.

If $Z$ is a~minimal center in $\mathbb{LCS}(X,D_{X})$ and $(X,
D_{X})$ is log canonical, then it follows from
Lemma~\ref{lemma:centers} that
$$
Z\cap g\big(Z\big)\ne \varnothing\iff Z=g\big(Z\big)%
$$
for every $g\in\bar{G}$.

\begin{lemma}
\label{lemma:Kawamata-Shokurov-trick} Suppose that $Z$ is
a~minimal center in $\mathbb{LCS}(X, D_{X})$, the log pair $(X,
D_{X})$ is log canonical, and $D_{X}$ is ample. Let $\epsilon$ be
an~arbitrary rational number such that $\epsilon>1$. Then there
exists an~effective $\bar{G}$-in\-va\-riant
$\mathbb{Q}$-divisor~$D$ on~the~variety~$X$ such~that
$$
\mathbb{LCS}\Big(X, D\Big)=\bigcup_{g\in\bar{G}}\Big\{g\big(Z\big)\Big\},%
$$
the~log pair $(X,D)$ is log canonical, and the~equivalence
$D\sim_{\mathbb{Q}} \epsilon D_{X}$ holds.
\end{lemma}

\begin{proof}
Take $m\in\mathbb{Z}$ such that $mD_{X}$ is a~very ample Cartier
divisor. Take a~general~divisor $R$ in the linear system
$|nmD_{X}|$ such that $Z\subset\mathrm{Supp}(R)$ and $R$ is
$\bar{G}$-invariant, where $n\gg 0$. Then
$$
\bigcup_{g\in\bar{G}}\Big\{g\big(Z\big)\Big\}\subseteq\mathbb{LCS}\Big(X, \lambda D_{X}+\mu R\Big)\subseteq\mathbb{LCS}\big(X, D_{X}\big)%
$$
for some positive rational numbers $\lambda$ and $\mu$ such that
$\lambda<1\leqslant \lambda+\mu nm<\epsilon$. One has $\lambda
D_{X}+\mu R\sim_{\mathbb{Q}} (\lambda+\mu nm)D_{X}$.

It follows from the~generality of the~divisor $R$ that $(X, \mu
R)$ is Kawamata log terminal, and
$$
\mathbb{LCS}\Big(X, \lambda D_{X}+\mu R\Big)=\bigcup_{g\in\bar{G}}g\big(Z\big),%
$$
because $\lambda<1$ and $n\gg 0$. Then there is
$\theta\in\mathbb{Q}_{>0}$ such that
$0<1-\theta\mu\leqslant\lambda<1$ and
$$
\bigcup_{g\in\bar{G}}\Big\{g\big(Z\big)\Big\}\subseteq\mathbb{LCS}\Big(X, \big(1-\theta\mu\big)D_{X}+\mu R\Big)\subseteq\mathbb{LCS}\Big(X, \lambda D_{X}+\mu R\Big),%
$$
but the~log pair $(X,(1-\theta\mu)D_{X}+\mu R)$ is log canonical
at the~general point of $Z$.

Note that for a~fixed $R$, the~number $\theta$ is a~function of
$\mu$. In the~above process, we can choose the~number $\mu$ so
that $1\leqslant 1-\theta\mu+\mu nm<\epsilon$ and
$$
\mathbb{LCS}\Big(X, \big(1-\theta\mu\big)D_{X}+\mu R\Big)=\bigcup_{g\in\bar{G}}\Big\{g\big(Z\big)\Big\},%
$$
because $Z$ is a~minimal center in $\mathbb{LCS}(X,D_{X})$ (see
Lemma~\ref{lemma:centers}). Put
$$
D=\big(1-\theta\mu\big)D_{X}+\mu R+\frac{\epsilon-1-\theta\mu+\mu nm}{nm}M,%
$$
where $M$ is a~general $\bar{G}$-invariant divisor in $|R|$. Then
$D$ is the~required divisor.
\end{proof}

Suppose now that $X=\mathbb{P}^{n}$. In this case we can say much
more about the locus $\mathrm{LCS}(X,D_{\mathbb{P}^{n}})$ and the
set $\mathbb{LCS}(\mathbb{P}^{n},D_{\mathbb{P}^{n}})$.

\begin{lemma}
\label{lemma:degree}  Let $H$ be a~hyperplane in $\mathbb{P}^{n}$,
and let $\mu$ be a~non-negative rational number such that
$D_{\mathbb{P}^{n}}\qlin \mu H$. Suppose that the locus
$\mathrm{LCS}(\mathbb{P}^{n},D_{\mathbb{P}^{n}})$
is~an~equidimensional subvariety in $\mathbb{P}^{n}$ of
codimension $s$. Put
$$
r=\left\{\aligned
&\lceil\mu-s-1\rceil\ \text{if}\ \mu\not\in\mathbb{Z},\\
&\lceil\mu-s-1\rceil+1\ \text {if}\ \mu\in\mathbb{Z}.\\
\endaligned
\right.
$$
Then $r\geqslant 0$ and
$$
\mathrm{deg}\Big(\mathrm{LCS}\big(\mathbb{P}^{n},D_{\mathbb{P}^{n}}\big)\Big)\leqslant {s+r\choose r}.%
$$
\end{lemma}

\begin{proof}
Put $Y=\mathrm{LCS}(\mathbb{P}^{n},D_{\mathbb{P}^{n}})$. Let
$\Pi\subset\mathbb{P}^{n}$ be a~general linear subspace of
dimension $s$. Put $D=D_{\mathbb{P}^{n}}\vert_{\Pi}$ and
$\Lambda=H\cap \Pi$. Then $\mathrm{deg}(Y)=|Y\cap \Pi|$ and
$\mathrm{LCS}(\Pi,D)=Y\cap\Pi$ by
Remark~\ref{remark:hyperplane-reduction}. One has $K_{\Pi}+D\qlin
(\mu-s-1)\Lambda$.

It follows from Theorem~\ref{theorem:Shokurov-vanishing} that
there is an exact sequence of cohomology groups
$$
0\longrightarrow H^{0}\Bigg(\mathcal{O}_{\Pi}\big(r\Lambda\big)\otimes\mathcal{I}\Big(\Pi,D\Big)\Bigg)\longrightarrow H^{0}\Big(\mathcal{O}_{\Pi}\big(r\Lambda\big)\Big)\longrightarrow H^{0}\Big(\mathcal{O}_{\mathcal{L}(\Pi,D)}\Big)\longrightarrow 0,%
$$
and $\mathrm{Supp}(\mathcal{L}(\Pi,D))=\mathrm{LCS}(\Pi,D)=Y\cap
\Pi\ne\varnothing$. Therefore, we see that $r\geqslant 0$ and
$$
\mathrm{deg}\big(Y\big)=\big|Y\cap \Pi\big|\leqslant h^{0}\Big(\mathcal{O}_{\mathcal{L}(\Pi,D)}\Big)\leqslant h^{0}\Big(\mathcal{O}_{\Pi}\big(r\Lambda\big)\Big)=h^{0}\Big(\mathcal{O}_{\mathbb{P}^{s}}\big(r\big)\Big)={s+r\choose r},%
$$
which completes the~proof.
\end{proof}

Let
$\phi\colon\GL_{n+1}(\mathbb{C})\to\mathrm{Aut}(\mathbb{P}^{n})\cong
\PGL_{n+1}(\mathbb{C})$ be the natural projection, and let $G$ be
a~finite subgroup in $\GL_{n+1}(\mathbb{C})$ such that
$\bar{G}=\phi(G)$.

\begin{remark}
\label{remark:semiinvariants-and-points} If $G$ does not have
semi-invariants of degree at most $k$, then every $\bar{G}$-orbits
in $\mathbb{P}^{n}$ contains at least $k+1$ points, because every
$\bar{G}$-orbit consisting of $s$ points defines a
$\bar{G}$-invariant hypersurface in $\mathbb{P}^{n}$ that is a
union of $s$ hyperplanes.
\end{remark}

\begin{lemma}
\label{lemma:semiinvariants-and-codimension-one} Let $H$ be
a~hyperplane in $\mathbb{P}^{n}$, and let $\mu$ be a~non-negative
rational number such that $D_{\mathbb{P}^{n}}\qlin \mu H$. Suppose
that $G$ does not have semi-invariants of degree at most $\lfloor
\mu\rfloor$. Then
$\mathbb{LCS}(\mathbb{P}^{n},D_{\mathbb{P}^{n}})$ does not contain
subvarieties in $\mathbb{P}^{n}$ of codimension $1$. If in
addition $\lfloor \mu\rfloor\leqslant n+1$ and the log pair
$(\mathbb{P}^{n},D_{\mathbb{P}^{n}})$ is log canonical, then
$\mathbb{LCS}(\mathbb{P}^{n},D_{\mathbb{P}^{n}})$ does not contain
points.
\end{lemma}

\begin{proof}
Suppose that $\mathbb{LCS}(\mathbb{P}^{n},D_{\mathbb{P}^{n}})$
contains an irreducible subvariety $Y\subset\mathbb{P}^{n}$ of
codimension $1$. Let $R$ be the $\bar{G}$-orbit of the subvariety
$Y$. Then
$$
D_{\mathbb{P}^{n}}=a R+\Delta
$$
for some rational number $a\geqslant 1$ and some effective
$\mathbb{Q}$-divisor $\Delta$ on $\mathbb{P}^{n}$. Since
\mbox{$D_{\mathbb{P}^{n}}\qlin \mu H$}, we see that $R$ is a
hypersurface in $\mathbb{P}^{n}$ of degree at most $\lfloor
\mu/a\rfloor\leqslant \lfloor \mu\rfloor$, which is impossible,
because $G$ does not have semi-invariants of degree at most
$\lfloor \mu\rfloor$.

We see that $\mathbb{LCS}(\mathbb{P}^{n},D_{\mathbb{P}^{n}})$ does
not contain subvarieties in $\mathbb{P}^{n}$ of codimension $1$.
Let us show that $\mathbb{LCS}(\mathbb{P}^{n},D_{\mathbb{P}^{n}})$
does not contain points provided that $\lfloor \mu\rfloor\leqslant
n+1$ and the log pair $(\mathbb{P}^{n},D_{\mathbb{P}^{n}})$ is log
canonical.

Suppose that  $\lfloor \mu\rfloor\leqslant n+1$, the log pair
$(\mathbb{P}^{n},D_{\mathbb{P}^{n}})$ is log canonical, and
$\mathbb{LCS}(\mathbb{P}^{n},D_{\mathbb{P}^{n}})$ contains a point
$P\in \mathbb{P}^{n}$. Let us show that these assumptions lead to
a contradiction.

Let $\Sigma$ be the $\bar{G}$-orbit of the point $P$, and let
$\epsilon$ be a rational number such that $\epsilon>1$ and
$\lfloor \epsilon\mu\rfloor\leqslant n+1$. Then it follows from
Lemma~\ref{lemma:Kawamata-Shokurov-trick} that there is
an~effective $\bar{G}$-invariant $\mathbb{Q}$-divisor $D$ on
$\mathbb{P}^{n}$~such~that $D\sim_{\mathbb{Q}}\epsilon\mu H$, the
log pair $(\mathbb{P}^{n}, D)$ is log canonical and
$\Sigma=\mathrm{LCS}(\mathbb{P}^{n},  D)$.

Since $\lfloor \epsilon\mu\rfloor\leqslant n+1$, it follows from
Theorem~\ref{theorem:Shokurov-vanishing} that
$$
H^{0}\Big(\mathcal{O}_{\mathbb{P}^{n}}\big(1\big)\otimes\mathcal{I}\big(\mathbb{P}^{n},D\big)\Big)=0,
$$
because $K_{\mathbb{P}^{n}}+D\sim_{\mathbb{Q}}(\epsilon\mu-n-1)H$
and $\epsilon\mu-n-1<1$. Therefore, it follows from the~exact
sequence of cohomology groups
$$
0\longrightarrow
H^{0}\Big(\mathcal{O}_{\mathbb{P}^{n}}\big(1\big)\otimes\mathcal{I}\big(\mathbb{P}^{n},D\big)\Big)\longrightarrow H^{0}\Big(\mathcal{O}_{\mathbb{P}^{n}}\big(1\big)\Big)\longrightarrow H^{0}\big(\mathcal{O}_{\Sigma}\big)\longrightarrow 0%
$$
that $|\Sigma|\leqslant n+1$, which is impossible because $G$ does
not have semi-invariants of degree at most $\lfloor
\mu\rfloor\leqslant n+1$.
\end{proof}

\begin{remark}
\label{remark:real-representation}  If $G$ is conjugate to
a~subgroup in $\GL_{n+1}(\mathbb{R})$, then the~subgroup $G$ has
an~invariant~of~degree~$2$, which implies that
$\mathrm{lct}(\mathbb{P}^{n},\bar{G})\leqslant 2/(n+1)$.
\end{remark}

\begin{remark}
\label{remark:stupid-action} If $Z$ is a $\bar{G}$-invariant, then
there is a~homomorphism $\xi\colon\bar{G}\to\mathrm{Aut}(Z)$ that
must be a~monomorphism provided that $Z$ is not contained in
a~linear subspace of $\mathbb{P}^{n}$, because eigenvectors that
correspond to a~fixed~eigenvalue of any matrix in
$\GL_{n+1}(\mathbb{C})$ form a~vector subspace in
$\mathbb{C}^{n+1}$.
\end{remark}

\begin{theorem}
\label{theorem:Hurwitz-bound} Let $C$ be a~smooth irreducible
curve of genus $g\geqslant 2$. Then $|\mathrm{Aut}(C)|\leqslant
84(g-1)$.
\end{theorem}

\begin{proof}
The required inequality is the~famous Hurwitz bound
(see~\cite[Theorem~3.17]{Bre00}).
\end{proof}

\section{Exceptionality criterion}
\label{section:plt-blow-up}

Let $X$ be a~variety, let $B_{X}$ be an~effective
$\mathbb{Q}$-divisor on $X$~such~that the~log pair $(X,B_{X})$~has
at most Kawamata log terminal singularities, and the~divisor
\mbox{$-(K_{X}+B_{X})$} is ample. Recall that $(X,B_{X})$ is usually
called a~\emph{log Fano variety}. Let~\mbox{$\bar{G}\subset\mathrm{Aut}(X)$}
be a~finite subgroup such that
the~divisor $B_{X}$ is $\bar{G}$-invariant.

\begin{definition}
\label{definition:G-threshold} The global $\bar{G}$-invariant log
canonical threshold of the~log Fano variety $(X,B_{X})$ is a real
number $\mathrm{lct}(X, B_{X},\bar{G})$ that can be defined as
$$
\mathrm{inf}\left\{\mathrm{c}\big(X,B_{X},D_{X}\big)\in\mathbb{Q}\ \left|\ %
\aligned &D_{X}\ \text{is a~$\bar{G}$-invariant $\mathbb{Q}$-Cartier effective $\mathbb{Q}$-divisor}\\
&\text{on the~variety $X$ such that}\ D_{X}\sim_{\mathbb{Q}} -\Big(K_{X}+B_{X}\Big)\\
\endaligned\right.\right\}.%
$$
\end{definition}

For simplicity, we  put
$\mathrm{lct}(X,B_{X},\bar{G})=\mathrm{lct}(X,\bar{G})$ if
$B_{X}=0$. Similarly, we put
$\mathrm{lct}(X,B_{X},\bar{G})=\mathrm{lct}(X,B_{X})$ if $\bar{G}$
is trivial. Finally, we put
$\mathrm{lct}(X,B_{X},\bar{G})=\mathrm{lct}(X)$ if $B_{X}=0$ and
$\bar{G}$ is trivial. Then it follows from
\cite[Theorem~A.3]{ChSh08c} that
$\mathrm{lct}(X,\bar{G})=\alpha_{\bar{G}}(X)$ if $X$ is smooth and
$B_{X}=0$ (see Definition~\ref{definition:alpha-invariant}).

\begin{remark}
\label{remark:quotient-lct} Suppose that $B_{X}=0$. Put
$V=X\slash\bar{G}$. Let $\theta\colon X\to V$ be the~quotient map.
Then
$$
K_{X}\sim_{\mathbb{Q}}\theta^{*}\Big(K_{V}+R_{V}\Big),
$$
where $R_{V}$ is a~ramification $\mathbb{Q}$-divisor of
the~morphism $\theta$. Note that $-(K_{V}+R_{V})$ is an~ample
$\mathbb{Q}$-Cartier divisor, and $(V,R_{V})$ is Kawamata log
terminal by \cite[Proposition~3.16]{Ko97}. Moreover, it follows
from \cite[Proposition~3.16]{Ko97} that
$\mathrm{lct}(X,\bar{G})=\mathrm{lct}(V,R_{V})$.
\end{remark}

\begin{example}
\label{example:lct-line} Suppose that $X\cong\mathbb{P}^{1}$. Then
$B_{X}=\sum_{i=1}^{n}a_{i}P_{i}$, where $P_{i}$ is a~point, and
$a_{i}\in\mathbb{Q}$ such that $0\leqslant a_{i}<1$.  We may
assume that $a_{0}\leqslant \ldots\leqslant a_{n}$. Then
$$
\mathrm{lct}\Big(X,B_{X}\Big)=\frac{1-a_{n}}{2-\sum_{i=1}^{n}a_{i}},
$$
where $\sum_{i=1}^{n}a_{i}<2$, because the~divisor
$-(K_{X}+B_{X})$ is ample. Moreover, it follows from
Remark~\ref{remark:quotient-lct} that
$\mathrm{lct}(X,\bar{G})=2/\lambda$, where $\lambda$ is the length
of a~$\bar{G}$-orbit of the~smallest length (cf.
Theorem~\ref{theorem:Shokurov-n-2}).
\end{example}

\begin{lemma}
\label{lemma:fixed-or-mobile} The global log canonical threshold
$\mathrm{lct}(X, B_{X},\bar{G})$ is equal to
$$
\mathrm{inf}\left\{\mathrm{c}\Big(X,B_{X},\sum_{i=1}^{r}a_{i}\mathcal{D}_{i}\Big)\ \left|\ %
\aligned &\text{$\mathcal{D}_{i}$ is a~linear system and $a_{i}\in \mathbb{Q}_{\geqslant 0}$}\\
&\text{for every $i\in\{1,\ldots,r\}$,
$\sum_{i=1}^{r}a_{i}\mathcal{D}_{i}$ is $\bar{G}$-invariant,}\\
&\text{and
$\sum_{i=1}^{r}a_{i}\mathcal{D}_{i}\sim_{\mathbb{Q}}-(K_{X}+B_{X})$}\\%
\endaligned\right.\right\}.%
$$
\end{lemma}

\begin{proof}
The required assertion follows from Definition~\ref{definition:G-threshold} and
\cite[Theorem~4.8]{Ko97}.
\end{proof}

In general, it is unknown whether $\mathrm{lct}(X,B_{X},\bar{G})$
is a rational number or not (cf. \cite[Question~1]{Ti90b}). Of
course, we expect that $\mathrm{lct}(X,B_{X},\bar{G})$ is
rational. Moreover, we expect the following to be true.

\begin{conjecture}
\label{conjecture:stabilization} There is an effective
$\bar{G}$-invariant $\mathbb{Q}$-divisor $D_{X}$ on $X$ such that
$\mathrm{lct}(X,B_{X},\bar{G})=\mathrm{c}(X,B_{X},D_{X})\in\mathbb{Q}$
and $D_{X}\sim_{\mathbb{Q}} -(K_{X}+B_{X})$.
\end{conjecture}

Let $(V\ni O)$ be a~germ of a~Kawamata log terminal singularity,
and let \mbox{$\pi\colon W\to V$} be~a~birational morphism such that
the~exceptional locus of $\pi$ consists of one irreducible divisor
$E\subset W$ such that $O\in\pi(E)$, the~log pair $(W,E)$ has
purely log terminal singularities (see
\cite[Definition~3.5]{Ko97}), and $-E$ is a~$\pi$-ample
$\mathbb{Q}$-Cartier divisor.

\begin{theorem}
\label{theorem:plt-blow-up} The birational morphism $\pi\colon
W\to V$ does exist.
\end{theorem}

\begin{proof}
Modulo the Log Minimal Model Program in dimension
$\mathrm{dim}(V)$, the existence of the morphism $\pi$ follows
from \cite[Proposition~2.9]{Pr98plt} in the case when $V$ has
$\mathbb{Q}$-factorial singularities. It follows from
\cite[Theorem~1.5]{Kud01} that the $\mathbb{Q}$-factoriality
condition in \cite[Proposition~2.9]{Pr98plt} can be removed.
Moreover, the proofs of \cite[Proposition~2.9]{Pr98plt} and
\cite[Theorem~1.5]{Kud01} only need the Log Minimal Model Program
for log pairs with big boundaries, which is proved now in
\cite{BCHM06}.
\end{proof}

We say that $\pi\colon W\to V$ is a~\emph{plt blow up} of
the~singularity $(V\ni O)$.

\begin{definition}[{\cite[Definition~4.1]{Pr98plt}}]
\label{definition:weakly-exceptional} We say that $(V\ni O)$ is
\emph{weakly-exceptional} if it has unique plt blow up.
\end{definition}

Weakly-exceptional Kawamata log terminal singularities do exist
(see \cite[Example~2.2]{Kud01}).

\begin{lemma}[{\cite[Corollary~1.7]{Kud01}}]
\label{lemma:weakly-exceptional} If  $(V\ni O)$ is
weakly-exceptional, then $\pi(E)=O$.
\end{lemma}

Let $R_{1},\ldots,R_{s}$ be irreducible components of
$\mathrm{Sing}(W)$ such that
$\mathrm{dim}(R_{i})=\mathrm{dim}(W)-2$ and $R_{i}\subset E$ for
every $i\in\{1,\ldots,s\}$. Put
$$
\mathrm{Diff}_{E}\big(0\big)=\sum_{i=1}^{s}\frac{m_{i}-1}{m_{i}}R_{i},
$$
where $m_{i}$ is the~smallest positive integer such that $m_{i}E$ is~Cartier at
a~general point of $R_{i}$.

\begin{lemma}[{\cite[Theorem~7.5]{Ko97}}]
\label{lemma:plt-log-Fano} The~variety $E$ is normal, and
$(E,\mathrm{Diff}_{E}(0))$ is Kawamata log terminal.
\end{lemma}

Therefore, if $\pi(E)=O$, then the~log pair
$(E,\mathrm{Diff}_{E}(0))$ is a~log Fano variety, because $-E$ is
$\pi$-ample.

\begin{theorem}[{\cite[Theorem~2.1]{Kud01}}]
\label{theorem:weakly-exceptional-criterion} The~singularity
$(V\ni O)$ is weakly-exceptional if~and~only~if $\pi(E)=O$ and
$\mathrm{lct}(E,\mathrm{Diff}_{E}(0))\geqslant 1$.
\end{theorem}

\begin{theorem}[{\cite[Theorem~4.9]{Pr98plt}}]
\label{theorem:exceptional-criterion} The~singularity $(V\ni O)$
is exceptional if and only~if $\pi(E)=O$ and
$\mathrm{c}(E,\mathrm{Diff}_{E}(0),D_{E})>1$ for every effective
$\mathbb{Q}$-divisor $D_{E}$ on the~variety $E$ such that
$D_{E}\sim_{\mathbb{Q}}-(K_{E}+\mathrm{Diff}_{E}(0))$.
\end{theorem}

In particular, we see that if the~assertion of
Conjecture~\ref{conjecture:stabilization} is true, then $(V\ni O)$
is exceptional if and only if $\pi(E)=O$ and
$\mathrm{lct}(E,\mathrm{Diff}_{E}(0))>1$~holds.

\begin{corollary}
\label{corollary:exceptional-weakly-exceptional} If $(V\ni O)$ is
exceptional, then $(V\ni O)$ is weakly-exceptional.
\end{corollary}

It should be pointed out that
Theorem~\ref{theorem:exceptional-criterion} is an applicable
criterion. For instance, it can be used to construct
exceptional singularities of any dimension.

\begin{example}
\label{example:Kollar-Boyer} Suppose that $(V\ni O)$ is
a~Brieskorn--Pham hypersurface singularity
$$
\sum_{i=0}^{n}x_{i}^{a_{i}}=0\subset\mathbb{C}^{n+1}\cong\mathrm{Spec}\Big(\mathbb{C}\big[x_{0},x_{1},\ldots,x_{n}\big]\Big),
$$
where $n\geqslant 3$ and $2\leqslant a_{0}<a_{1}<\ldots<a_{n}$.
Arguing as in the~proof of~\cite[Theorem~34]{KGB}, we~see~that it
follows from Theorem~\ref{theorem:exceptional-criterion} that
the~singularity $(V\ni O)$ is exceptional if
$$
1<\sum\limits_{i=0}^{n}\frac{1}{a_i}<1+\mathrm{min}\Bigg\{\frac{1}{a_{0}},\frac{1}{a_{1}},\ldots,\frac{1}{a_{n}}\Bigg\}
$$
and $a_{0},a_{1},\ldots,a_{n}$ are pairwise coprime. This is
 satisfied if $a_{0},a_{1},\ldots,a_{n}$
are primes and
\begin{equation}
\label{equation:restricted-KGB-condition}
\frac{1}{a_0}+\frac{1}{a_1}+\ldots+\frac{1}{a_{n-1}}<1<\frac{1}{a_0}+\frac{1}{a_1}+\ldots+\frac{1}{a_{n-1}}+\frac{1}{a_{n}}.%
\end{equation}

Let us use induction to construct the~$(n+1)$-tuple
$(a_{0},a_{1}\ldots,a_{n})$ such that $a_{0},a_{1},\ldots,a_{n}$
are prime integers, and the~$(n+1)$-tuple
$(a_{0},a_{1}\ldots,a_{n})$ satisfies
the~inequalities~\ref{equation:restricted-KGB-condition}.

If $n=3$, then the~four-tuple $(a_0, a_1, a_2, a_3)=(2, 3, 7, 41)$
satisfies
the~inequalities~\ref{equation:restricted-KGB-condition}.

Suppose that $n\geqslant 4$, and there are prime numbers
\mbox{$2\leqslant c_{0}<c_{1}<c_{2}<\ldots<c_{n-1}$} such~that
$$
\frac{1}{c_0}+\frac{1}{c_1}+\ldots+\frac{1}{c_{n-2}}<1<\frac{1}{c_0}+\frac{1}{c_1}+\ldots+\frac{1}{c_{n-2}}+\frac{1}{c_{n-1}},%
$$
and assume that $c_{n-1}>8$ is the~largest prime with these
properties (for the~fixed numbers $c_0, \ldots, c_{n-2}$). It follows from
$c_{n-1}>8$ that there are prime numbers $p_1, p_2$ and~$p_3$ such
that $c_{n-1}<p_1<p_2<p_3<2c_{n-1}$ (see
\cite[p.~209,~(18)]{Ram19}). Put $(a_{0},a_{1}\ldots,a_{n})=(c_0,
\ldots, c_{n-2}, p_2, p_3)$. Then
$$
\sum_{i=0}^{n-2}\frac{1}{a_i}+\frac{1}{p_2}<
\sum_{i=0}^{n-2}\frac{1}{a_i}+\frac{1}{p_1}\leqslant
1<\sum_{i=0}^{n-2}\frac{1}{c_i}+\frac{1}{2c_{n-1}}+\frac{1}{2c_{n-1}}
<\sum_{i=0}^{n-2}\frac{1}{a_i}+\frac{1}{p_2}+\frac{1}{p_3}%
$$
by the maximality assumption imposed on $c_{n-1}$. Hence
the~$(n+1)$-tuple $(a_{0},a_{1}\ldots,a_{n})$ satisfies
the~inequalities~\ref{equation:restricted-KGB-condition}, which
completes the construction\footnote{ Alternatively, one can use
the Sylvester sequence to construct $(a_0, \ldots, a_n)$
explicitly (suggested by S.\,Galkin).}.
\end{example}

\medskip

Suppose, in addition, that $(V\ni O)$ is a~quotient singularity
$\mathbb{C}^{n+1}\slash G$, where $n\ge 1$ and $G$ is
a~finite~subgroup in $\GL_{n+1}(\mathbb{C})$. Put
$\bar{G}=\phi(G)$, where
$\phi\colon\GL_{n+1}(\mathbb{C})\to\mathrm{Aut}(\mathbb{P}^{n})\cong
\PGL_{n+1}(\mathbb{C})$
is the~natural projection.

\begin{remark}
\label{remark:plt-blow-up-quotient} Let
$\eta\colon\mathbb{C}^{n+1}\to V$ be the~quotient map. Then there
is a~commutative diagram
$$
\xymatrix{
U\ar@{->}[d]_{\gamma}\ar@{->}[rr]^{\omega}&&W\ar@{->}[d]^{\pi}\\
\mathbb{C}^{n+1}\ar@{->}[rr]_{\eta}&&V,}
$$
where $\gamma$ is the~blow up of $O$, the~morphism $\omega$ is
the~quotient map that is induced by the~lifted action~of $G$ on
the~variety $U$, and $\pi$ is a~birational morphism. Moreover, $\pi$ is
a~plt blow~up of the~singularity~$\C^{n+1}\slash G$.
\end{remark}

Thus, to prove the~existence of a~plt blow up of the quotient
singularity $\C^{n+1}\slash G$ we do not need to use
Theorem~\ref{theorem:plt-blow-up}.

\begin{theorem}
\label{theorem:weakly-exceptional-quotient}
Suppose that the group $G\subset\GL_{n+1}(\C)$
does not contain reflections.
Then the~singularity $\C^{n+1}\slash G$ is weakly-exceptional $\iff$
$\mathrm{lct}(\mathbb{P}^{n},\bar{G})\geqslant 1$.
\end{theorem}

\begin{proof}
Let us use the~notation and assumptions of
Remark~\ref{remark:plt-blow-up-quotient}. Let $F$ be
the~exceptional divisor of the~blow up $\gamma$. Put
$E=\omega(F)$. Then $F\cong\mathbb{P}^{n}$ and
$E\cong\mathbb{P}^{n}\slash\bar{G}$. Since the group $G$ does not
contain reflections, it follows from
Remark~\ref{remark:quotient-lct} that
$\mathrm{lct}(\mathbb{P}^{n},\bar{G})=\mathrm{lct}(E,\mathrm{Diff}_{E}(0))$,
which implies that the singularity~$\C^{n+1}\slash G$
is weakly-exceptional if and only if
$\mathrm{lct}(\mathbb{P}^{n},\bar{G})\geqslant 1$ by
Theorem~\ref{theorem:exceptional-criterion}.
\end{proof}

\begin{theorem}
\label{theorem:exceptional-quotient} Suppose that the group
$G\subset\GL_{n+1}(\C)$ does not contain reflections.
Then the singularity~$\C^{n+1}\slash G$
is exceptional if
and only if for any $\bar{G}$-invariant effective
$\mathbb{Q}$-divisor $D$ on $\mathbb{P}^{n}$ such that
$D\sim_{\mathbb{Q}} -K_{\mathbb{P}^{n}}$
the~log pair $(\mathbb{P}^{n}, D)$ is Kawamata log terminal.%
\end{theorem}

\begin{proof}
Arguing as in the~proof of
Theorem~\ref{theorem:weakly-exceptional-quotient} and using
Theorem~\ref{theorem:exceptional-criterion} together with
\cite[Proposition~3.16]{Ko97}, we obtain the required assertion.
\end{proof}

Recall that the subgroup $G\subset\GL_{n+1}(\mathbb{C})$ is
said to be transitive if the~corresponding $(n+1)$-dimensional
representation is irreducible (see {\cite{Bli17}). Note that $G$
is transitive if it is primitive. As an easy application of
Theorems~\ref{theorem:exceptional-quotient} and
\ref{theorem:weakly-exceptional-quotient} in conjunction with
Lemma~\ref{lemma:fixed-or-mobile} one can establish the relation
between the primitivity of the group~$G$ (transitivity, respectively)
and the exceptionality of the singularity~$\C^{n+1}\slash G$
(weak-exceptionality,
respectively).

\begin{theorem}\label{theorem:primitive}
Suppose that the group $G\subset\GL_{n+1}(\C)$
is not primitive (not transitive, respectively).
Then there exists a $\bar{G}$-invariant effective $\Q$-divisor $D$
on $\P^n$ such that $D\sim_{\mathbb{Q}} -K_{\P^n}$ and
the pair $(\P^n,D)$ is not Kawamata log terminal (not log canonical,
respectively).
\end{theorem}

\begin{proof}
We will only prove that if the group $G$
is not primitive, then
there exists a
$\bar{G}$-invariant effective $\Q$-divisor $D$
on $\P^n$ such that $D\sim_{\mathbb{Q}} -K_{\P^n}$ and
the pair $(\P^n,D)$ is not Kawamata log terminal, since
the remaining assertion can be proved similarly.

Suppose that $G$ is not primitive. Then there is~a~non-trivial
decomposition
$$
\mathrm{Spec}\Big(\mathbb{C}\big[x_{0},x_{1},\ldots,x_{n}\big]\Big)\cong\mathbb{C}^{n+1}=\bigoplus_{i=1}^{r}V_{i}
$$
such that $g(V_{i})=V_{j}$ for all $g\in G$. We may assume that
$\mathrm{dim}(V_{1})\leqslant\ldots\leqslant\mathrm{dim}(V_{r})$.

Put $d=\mathrm{dim}(V_{1})$. Then $d\leqslant\lfloor
(n+1)/2\rfloor$. We may assume that $V_{1}\subset\mathbb{C}^{n+1}$
is given by $x_{d}=x_{d+1}=x_{d+2}=\ldots=x_{n}=0$. Let
$\mathcal{M}_{1}$ be a~linear system on $\mathbb{P}^{n}$ that
consists of hyperplanes that are given by
$$
\sum_{i=0}^{d-1}\lambda_{i}x_{i}=0\subset\mathbb{P}^{n}\cong\mathrm{Proj}\Big(\mathbb{C}\big[x_{0},x_{1},\ldots,x_{n}\big]\Big),
$$
where $\lambda_{i}\in\mathbb{C}$. Let
$\mathcal{M}_{1},\ldots,\mathcal{M}_{s}$ be the~$\bar{G}$-orbit of
the~linear system
 $\mathcal{M}_{1}$. Then
$$
\frac{n+1}{s}\Bigg(\sum_{i=1}^{s}\mathcal{M}_{i}\Bigg)\sim_{\mathbb{Q}} -K_{\mathbb{P}^{n}},%
$$
where $s\leqslant\lfloor (n+1)/d\rfloor$. Let
$\Lambda\subset\mathbb{P}^{n}$ be a~linear subspace that is given
by the equations \mbox{$x_{0}=\ldots=x_{d}=0$}. Then
$$
\frac{n+1}{s}\mathrm{mult}_{\Lambda}\Bigg(\sum_{i=1}^{s}\mathcal{M}_{i}\Bigg)\geqslant
\frac{n+1}{s}\mathrm{mult}_{\Lambda}\big(\mathcal{M}_{1}\big)=\frac{n+1}{s}\geqslant
d=n-\mathrm{dim}\big(\Lambda\big),
$$
which implies the desired assertion by
Lemma~\ref{lemma:fixed-or-mobile}.
\end{proof}

\begin{corollary}
\label{corollary:primitive-baby} Suppose that the group
$G\subset\GL_{n+1}(\C)$ is not
primitive (not transitive, respectively). Then
$\mathrm{lct}(\P^n,\bar{G})\leqslant 1$
($\mathrm{lct}(\P^n,\bar{G})<1$, respectively).
\end{corollary}

Applying Theorems~\ref{theorem:weakly-exceptional-quotient},
\ref{theorem:exceptional-quotient}
and~\ref{theorem:primitive}, we obtain the following.

\begin{corollary}[{\cite[Proposition~2.1]{Pr00}}]
\label{corollary:primitive}
Suppose that the group $G\subset\GL_{n+1}(\C)$ does not contain reflections.
Then the group~$G$ is primitive
(transitive, respectively) provided that the singularity~$\C^{n+1}\slash G$
is exceptional
(weakly-exceptional, respectively).
\end{corollary}

Let us show how to apply
Theorems~\ref{theorem:weakly-exceptional-quotient} and
\ref{theorem:exceptional-quotient} (cf.
\cite[Example~1.9]{Ch07b}).

\begin{theorem}
\label{theorem:Dima-Yura-strong} Suppose that $G\subset\GL_{3}(\C)$.
Then  $\mathrm{lct}(\mathbb{P}^{2},\bar{G})\geqslant 4/3$ if and
only if $G$ does not have semi-invariants of degree at most $3$.
\end{theorem}

\begin{proof}
Suppose that the~subgroup $G$ does not have semi-invariants of
degree at most~$3$. To complete the~proof we must~show~that
$\mathrm{lct}(\mathbb{P}^{2},\bar{G})\geqslant 4/3$, because the
remaining implication is obvious.

Suppose that the~strict inequality
$\mathrm{lct}(\mathbb{P}^{2},\bar{G})<4/3$ holds. Then there exist
a~positive rational number
$\lambda<4/3$~and~an~effective~$\bar{G}$-invariant
$\mathbb{Q}$-divisor $D$ on $\mathbb{P}^{2}$ such that
$D\sim_{\mathbb{Q}}-K_{\mathbb{P}^{2}}$, and the log pair
$(\mathbb{P}^{2},\lambda D)$ is strictly log canonical. Applying
Lemma~\ref{lemma:semiinvariants-and-codimension-one}, we obtain a
contradiction.
\end{proof}

Using Theorems~\ref{theorem:exceptional-quotient}
and~\ref{theorem:Dima-Yura-strong}, we obtain the following.

\begin{corollary}
\label{corollary:Dima-Yura-strong}
Suppose that the group $G\subset\GL_3(\C)$ does not contain reflections.
Then the~following  are equivalent:
\begin{itemize}
\item the~singularity $\C^{3}\slash G$ is exceptional,%

\item the~subgroup $G$ does not have semi-invariants of degree at most $3$,%

\item the~inequality $\mathrm{lct}(\mathbb{P}^{2},\bar{G})\geqslant 4/3$ holds.%
\end{itemize}
\end{corollary}

Arguing as in the~proof of Theorem~\ref{theorem:Dima-Yura-strong},
we easily obtain a similar assertion that can be used for the
classification of three-dimensional weakly exceptional quotient
singularities (see \cite{Sak10}).

\begin{theorem}
\label{theorem:Yura-Dima-weakly-exceptional} Suppose that
$G\subset\GL_{3}(\C)$. Then
$\mathrm{lct}(\mathbb{P}^{2},\bar{G})\geqslant 1$ if and only if
$G$ does not have semi-invariants of degree at most $2$.
\end{theorem}

\begin{proof}
The proof is left to the reader.
\end{proof}

\medskip

Suppose that $n+1=2l$ for some integer $l\geqslant 2$. Let
$G_{1}\subset\SL_2(\mathbb{C})$ and
\mbox{$G_{2}\subset\SL_l(\mathbb{C})$} be finite subgroups, let
$\mathbb{M}$ be the~vector space
of~$2\times l$-matrices with entries in~$\mathbb{C}$.
For every $(g_{1},g_{2})\in G_{1}\times G_{2}$ and
every $M\in\mathbb{M}$, put
$$
\Big(g_{1},g_{2}\Big)\big(M\big)=g_{1}Mg_{2}^{-1}
\in\mathbb{M}\cong\mathbb{C}^{2l},%
$$
which induces a~homomorphism $\varphi\colon G_{1}\times
G_{2}\to\SL_{2l}(\mathbb{C})$. Note that
$|\mathrm{ker}(\varphi)|\leqslant 2$ if $n$ is even, and $\varphi$
is a~monomorphism if $n$ is odd.

\begin{lemma}
\label{lemma:Segre} Suppose that $G=\varphi(G_{1}\times G_{2})$.
Then $\mathrm{lct}(\mathbb{P}^{n}, \bar{G})<1$.
\end{lemma}

\begin{proof}
Put $s=l-1$. Let
$\psi\colon\mathbb{P}^{1}\times\mathbb{P}^{s}\to\mathbb{P}^{n}$ be
the~Segre embedding. Put
$Y=\psi(\mathbb{P}^{1}\times\mathbb{P}^{s})$ and let $\mathcal{Q}$
be the~linear system consisting of all quadric hypersurfaces in
$\mathbb{P}^{n}$ that pass through the~subvariety $Y$. Then
$\mathcal{Q}$ is a~non-empty $\bar{G}$-invariant linear system.
The log pair $(\mathbb{P}^{n}, l\mathcal{Q})$ is not log-canonical
along $Y$, which implies that $\mathrm{lct}(\mathbb{P}^{n},
\bar{G})<1$ by Lemma~\ref{lemma:fixed-or-mobile}.
\end{proof}

As an application of Lemma~\ref{lemma:Segre} one obtains non-exceptionality
of some quotient singularities.

\begin{example}[{cf. Theorem~\ref{theorem:Vanya-Kostya-invariants}}]
\label{example:Severi} Suppose that $G=\varphi(G_{1}\times G_{2})$
and $l=3$. Then
the singularity~$\C^{6}\slash G$ is not exceptional
by Theorem~\ref{theorem:criterion} and Lemma~\ref{lemma:Segre}.
On the other hand, if $G_{1}\cong 2.\A_5$ and $G_2\cong 3.\A_6$,
then $G$ has no semi-invariants
of degree at most $6$ which can be shown by direct computation.
\end{example}

Suppose that $l=2$. The~transposition of matrices in $\mathbb{M}$
induces an~involution $\iota\in\SL_4(\mathbb{C})$.

\begin{lemma}
\label{lemma:Segre-2} If $G$ is generated by $\varphi(G_{1}\times
G_{2})$ and $\iota$, then $\mathrm{lct}(\mathbb{P}^{3},
\bar{G})<1$.
\end{lemma}

\begin{proof}
See the~proof of Lemma~\ref{lemma:Segre}.
\end{proof}

\section{Four-dimensional case}
\label{section:4}

Shokurov (see~\cite{Sho93}) and Prokhorov and Markushevich (see~\cite{MarPr99})
obtained an explicit classification of
exceptional quotient singularities of dimension $2$ and $3$.
Namely, for Gorenstein quotient singularities they prove the following.

\begin{theorem}[{\cite[Example~5.2.3]{Sho93}}]
\label{theorem:Shokurov-n-2}  Let $G$ be the finite subgroup in
$\SL_2(\C)$. Then the singularity~\mbox{$\C^2/G$}
is exceptional if and only if~$G$
is a binary central extension of one of the following groups:
$\A_{4}$, $\SS_{4}$ or $\A_{5}$.
\end{theorem}

\begin{theorem}[{\cite[Theorem~3.13]{MarPr99}}]
\label{theorem:Dima-Yura-n-3} Let $G$ be a finite subgroup in
$\SL_3(\C)$. Then the singularity~\mbox{$\C^3/G$}
is exceptional if and only if
$G$ is one of the following subgroups:
\begin{itemize}
\item a central extension of $\PSL_2(\mathbb{F}_{7})$, which is isomorphic to either $\PSL_2(\mathbb{F}_{7})$ or $\mathbb{Z}_{3}\times \PSL_2(\mathbb{F}_{7})$,%
\item a non-trivial central extension~$3.\A_6$ of the alternating
group~$\A_{6}$ by~$\mathbb{Z}_{3}$, %
\item the~Hessian group, which can be characterized by the~exact sequence%
$$
1\longrightarrow \mathbb{H}\Big(3,\mathbb{F}_{3}\Big)\longrightarrow G\longrightarrow\SS_{4}\longrightarrow 1,%
$$
where $\mathbb{H}(3,\mathbb{F}_{3})$ is the~Heisenberg group
consisting of all
unipotent $3\times 3$-matrices with entries in $\mathbb{F}_{3}$,%
\item the~normal subgroup of the~Hessian group of index $3$ that contains $\mathbb{H}(3, \mathbb{F}_{3})$.%
\end{itemize}
\end{theorem}

The purpose of this section is to present an analogous
classification for exceptional singularities of dimension $4$
(see Theorem~\ref{theorem:Vanya-Kostya-n-4}), and prove some relevant results.

\medskip

Let $\bar{G}$ be a~finite subgroup in
$\mathrm{Aut}(\mathbb{P}^{3})$, and let
$\phi\colon\GL_4(\mathbb{C})\to\mathrm{Aut}(\mathbb{P}^{3})$
be the~natural projection. Then there is a~finite subgroup in
$\SL_4(\mathbb{C})$ such that $\phi(G)=\bar{G}$. Moreover,
if $G$ is primitive, then it follows from \cite{Bli17} and
\cite{Fe71} that one may assume that $Z(G)\subseteq[G,G]$, where
$Z(G)$ and $[G,G]$ are the~center and the~commutator of the group
$G$, respectively.

As a warming-up we start with a result that can be applied to a
classification of four-dimensional weakly exceptional quotient
singularities (see \cite{Sak10}).

\begin{theorem}
\label{theorem:Vanya-Kostya-SL-4-weakly-exceptional}
The~inequality $\mathrm{lct}(\mathbb{P}^{3},\bar{G})\geqslant 1$
holds  if and only if the~following three conditions are
satisfied: the~group $G$ is transitive, the~group $G$ does not
have semi-invariants of degree at most $3$, and\footnote{
One can show that the third condition of
Theorem~\ref{theorem:Vanya-Kostya-SL-4-weakly-exceptional} is not redundant.
Namely, if $G\subset\SL_4(\C)$ is a primitive group isomorphic to~$2.\A_5$,
then~$G$ has no semi-invariants of degree at most~$3$, but there is
a $\bar{G}$-invariant twisted cubic in~$\P^3$. In fact, the
primitive group $G\cong 2.\A_5$ gives essentially
the only example of this kind.}
there is no
$\bar{G}$-invariant smooth rational cubic curve in
$\mathbb{P}^{3}$.
\end{theorem}

\begin{proof}
Let us prove the~$\Rightarrow$-part. If $G$ has a~semi-invariant
of degree at most $3$, then
$\mathrm{lct}(\mathbb{P}^{3},\bar{G})\leqslant 3/4$ by
Definition~\ref{definition:G-threshold}. If $G$ is not transitive,
then $\mathrm{lct}(\mathbb{P}^{3},\bar{G})<1$ by
Corollary~\ref{corollary:primitive-baby}.

Suppose that there is a~$\bar{G}$-invariant smooth rational cubic
curve $C\subset\mathbb{P}^{3}$. Let
$R\subset\nolinebreak\mathbb{P}^{3}$ be~the~surface that is swept
out by lines that are tangent to $C$. Then
$\mathrm{c}(\mathbb{P}^{3}, R)=5/6$ the~surface $R$ is
$\bar{G}$-invariant, and $\mathrm{deg}(R)=4$. Hence, we see that
$\mathrm{lct}(\mathbb{P}^{3},\bar{G})\le 5/6$.

Let us prove the~$\Leftarrow$-part. Suppose that $G$ is
transitive, the~subgroup $G$ has no semi-invariants of
degree~at~most~$3$, there is no $\bar{G}$-invariant smooth
rational cubic curve in $\mathbb{P}^{3}$, but
$\mathrm{lct}(\mathbb{P}^{3},\bar{G})<1$.

There is an~effective $\bar{G}$-invariant $\mathbb{Q}$-divisor $D$
on $\mathbb{P}^{3}$ such that
$D\sim_{\mathbb{Q}}-K_{\mathbb{P}^{3}}$ and a~positive rational
number $\lambda<1$ such~that $(\mathbb{P}^{3},\lambda D)$ is
strictly log canonical. Let $S$ be an irreducible
subvariety of $\mathbb{P}^{3}$ that is a~minimal center in
$\mathbb{LCS}(\mathbb{P}^{3},\lambda D)$. By
Lemma~\ref{lemma:Kawamata-Shokurov-trick}, we may assume that
$$
\mathbb{LCS}\Big(\mathbb{P}^{3},\lambda D\Big)=\bigcup_{g\in\bar{G}}\Big\{g\big(S\big)\Big\},%
$$
where $\mathrm{dim}(S)\ne 2$, because $G$ has no semi-invariants
of degree at most~$3$.

The locus $\mathrm{LCS}(\mathbb{P}^{3},\lambda D)$ is connected by
Corollary~\ref{corollary:connectedness}. Then $S$ is
$\bar{G}$-invariant by Lemma~\ref{lemma:centers}. Since the group
$G$ is transitive, we see that $S$ is not a point. We see that $S$
is a~curve. Then $\mathrm{deg}(S)\leqslant 3$ by
Lemma~\ref{lemma:degree}, and $S$ is not contained in a~plane,
because $G$ is transitive. Hence $S$ is a~smooth rational cubic
curve.
\end{proof}

Combining Remark~\ref{remark:stupid-action},
Theorem~\ref{theorem:Vanya-Kostya-SL-4-weakly-exceptional} and
the~classification of finite subgroups in
$\PGL_2(\mathbb{C})$, we easily obtain the~following
result (cf. Theorem~\ref{theorem:Yura-Dima-weakly-exceptional}).

\begin{corollary}
\label{corollary:Vanya-Kostya-SL-4-weakly-exceptional}
The~inequality $\mathrm{lct}(\mathbb{P}^{3},\bar{G})\geqslant 1$
holds if the~following three conditions are satisfied: the~group $G$ is
transitive, the~group $G$ does not have semi-invariants of degree at most $3$,
and the group~$\bar{G}$ is not isomorphic to the alternating group~$\A_5$.
\end{corollary}

\medskip
The main purpose of this section is to prove the~following result
(cf. Theorem~\ref{theorem:Dima-Yura}).

\begin{theorem}
\label{theorem:Vanya-Kostya} The~inequality
$\mathrm{lct}(\mathbb{P}^{3},\bar{G})\geqslant 5/4$ holds if
the~following three conditions are satisfied: the~group $G$ is
primitive, the~group $G$ does not have semi-invariants of degree
at most $4$, and the~inequality $|\bar{G}|\geqslant 169$ holds.
\end{theorem}

\begin{proof}
Suppose that $G$ is primitive and does not have semi-invariants of
degree at most $4$, the~inequality $|\bar{G}|\geqslant 169$ holds,
but $\mathrm{lct}(\mathbb{P}^{3},\bar{G})<5/4$. Let us derive
a~contradiction.

There is an~effective $\bar{G}$-invariant $\mathbb{Q}$-divisor $D$
on $\mathbb{P}^{3}$~such~that
$D\sim_{\mathbb{Q}}-K_{\mathbb{P}^{3}}$ and a~positive rational
number $\lambda<5/4$ such~that $(\mathbb{P}^{3},\lambda D)$ is
strictly log canonical.

Let $S$ be an irreducible subvariety in $\mathbb{P}^{3}$ that is
a~minimal center in the~set $\mathbb{LCS}(\mathbb{P}^{3},\lambda
D)$. Then $S$ is a curve by
Lemma~\ref{lemma:semiinvariants-and-codimension-one}.

Note that $g(S)\in\mathbb{LCS}(\mathbb{P}^{3},\lambda D)$ for
every $g\in\bar{G}$, because the~divisor~$D$ is
$\bar{G}$-invariant. It follows from Lemma~\ref{lemma:centers}
that
$$
S\cap g\big(S\big)\ne \varnothing\iff S=g\big(S\big)%
$$
for every $g\in\bar{G}$. It follows from
Lemma~\ref{lemma:Kawamata-Shokurov-trick} that we may assume that
$$
\mathbb{LCS}\Big(\mathbb{P}^{3},\lambda D\Big)=\bigcup_{g\in\bar{G}}\Big\{g\big(S\big)\Big\}.%
$$

Let $\mathcal{I}$ be the~multiplier ideal sheaf of the~log pair
$(\mathbb{P}^{3},\lambda D)$, and let $\mathcal{L}$ be the~log
canonical singularities subscheme of the~log pair
$(\mathbb{P}^{3},\lambda D)$. Then there is an~exact sequence
\begin{equation}
\label{equation:exact-sequence-1}
0\longrightarrow
H^{0}\Big(\mathcal{O}_{\mathbb{P}^{3}}\big(1\big)\otimes\mathcal{I}\Big)\longrightarrow H^{0}\Big(\mathcal{O}_{\mathbb{P}^{3}}\big(1\big)\Big)\longrightarrow H^{0}\Big(\mathcal{O}_{\mathcal{L}}\otimes\mathcal{O}_{\mathbb{P}^{3}}\big(1\big)\Big)\longrightarrow 0%
\end{equation}
by Theorem~\ref{theorem:Shokurov-vanishing}. Then it follows from
Theorem~\ref{theorem:Kawamata} that $S$ is a~smooth curve of genus
$g$ such that $2g-2<\mathrm{deg}(S)$.

Let $Z$ be the~$\bar{G}$-orbit of the~curve $S$. Then $Z$ is
smooth and $\mathrm{deg}(Z)\leqslant 6$ by
Lemma~\ref{lemma:degree}.~Then $2g-2<\mathrm{deg}(S)\leqslant 6$,
which implies that $g\leqslant 3$. Note that $Z=\mathcal{L}$ by
Remark~\ref{remark:log-canonical-subscheme}, because
$(\mathbb{P}^{3},\lambda D)$ is log canonical. Moreover, the curve
$Z$ is not contained in a~plane, because $G$ is transitive.

Let $r$ be the~number of irreducible components of $Z$. Then
$6\geqslant \mathrm{deg}(Z)=r\mathrm{deg}(S)$, which implies that
$r\leqslant 6$. Note that $g=0$ if $r\geqslant 3$.

Using $(\ref{equation:exact-sequence-1})$ and the~Riemann--Roch
theorem, we see~that
\begin{equation}
\label{equation:inequalities-4} 4=h^{0}\Big(\mathcal{O}_{\mathcal{L}}\otimes\mathcal{O}_{\mathbb{P}^{3}}\big(1\big)\Big)=r\Big(\mathrm{deg}\big(S\big)-g+1\Big),%
\end{equation}
because $\mathcal{L}=Z$ and $2g-2<\mathrm{deg}(S)$. In particular,
we see that $r\leqslant 2$.

One has $\mathrm{deg}(S)\ne 1$, because $G$ is primitive. Thus $S$
is not contained in a~plane, because otherwise the~$\bar{G}$-orbit
of the~plane spanned by $S$ would give a~semi-invariant of $G$ of
degree~$1$~or~$2$. Thus, we have $6\geqslant
\mathrm{deg}(Z)=r\mathrm{deg}(S)\geqslant 3r$.

If $r=2$, then $\mathrm{deg}(S)=3$ and $g=0$, which contradicts
the~equalities~\ref{equation:inequalities-4}. We see that $r=1$
and $Z=S$. Then $g\leqslant 1$ by
Theorem~\ref{theorem:Hurwitz-bound} and
Remark~\ref{remark:stupid-action}, because $|\bar{G}|\geqslant
169$.

Arguing as in the~proof of
Theorem~\ref{theorem:Vanya-Kostya-SL-4-weakly-exceptional}, we see
that $g\ne 0$, because $G$ does not have semi-invariants of degree
$4$. Then it follows from $(\ref{equation:inequalities-4})$ that
$g=1$ and $\mathrm{deg}(S)=4$. We see that $S=Q_{1}\cap Q_{2}$,
where $Q_{1}$ and $Q_{2}$ are irreducible quadrics in~$\mathbb{P}^{3}$.

Let $\mathcal{P}$ be a~pencil generated~by $Q_{1}$~and~$Q_{2}$.
Then $\mathcal{P}$ contains exactly $4$ singular surfaces, which
are simple quadric cones. This means that there is a
$\bar{G}$-orbit in $\mathbb{P}^{3}$ consisting of at most $4$
points, which~is~impossible by
Remark~\ref{remark:semiinvariants-and-points}.
\end{proof}

In the rest of this section we will refine the assertion of
Theorem~\ref{theorem:Vanya-Kostya} by removing the assumption that
$\bar{G}$ contains at least $169$ elements and providing an
explicit list of possible finite subgroups in
$\PGL_{4}(\mathbb{C})$ that satisfy all hypothesis of
Theorem~\ref{theorem:Vanya-Kostya} (cf.
Theorems~\ref{theorem:Shokurov-n-2} and
\ref{theorem:Dima-Yura-n-3}). Let us start with the following
example.

\begin{example}[{see \cite[\S 123]{Bli17},
\cite{Nie92}}] \label{example:dim-4-primitive} Let $\mathbb{H}$ be
a~subgroup in $\SL_4(\mathbb{C})$ that is conjugate to the subgroup
generated by
$$
\left(\begin{array}{cccc}
0 & 0 & 1&0 \\
0 & 0 & 0&1 \\
1 & 0 & 0&0 \\
0 & 1 & 0&0 \\
\end{array}\right),
\left(\begin{array}{cccc}
0 & 1 & 0&0 \\
1 & 0 & 0&0 \\
0 & 0 & 0&1 \\
0 & 0 & 1&0 \\
\end{array}
\right), \left(\begin{array}{cccc}
1 & 0 & 0&0 \\
0 & 1 & 0&0 \\
0 & 0 & -1&0 \\
0 & 0 & 0&-1 \\
\end{array}
\right), \left(\begin{array}{cccc}
1 & 0 & 0&0 \\
0 & -1 & 0&0 \\
0 & 0 & 1&0 \\
0 & 0 & 0&-1 \\
\end{array}
\right),
$$
and let $N\subset\SL_4(\mathbb{C})$ be the~normalizer of
the~subgroup $\mathbb{H}$. There is an~exact sequence of
groups\footnote{The choice of the epimorphism $\beta$ is not
canonical even up to conjugation, due to the existence of outer
automorphisms of $\SS_6$. There are essentially two possible
choices of $\beta$. To fix one of them we use the fact that the
subspace $W\subset\mathrm{Sym}^4(\C^4)$ of
$\tilde{\mathbb{H}}$-invariant quartics is five-dimensional;
moreover, the group $N\slash\tilde{\mathbb{H}}$ acts on $W$, and
$W$ is an irreducible representation of
$N\slash\tilde{\mathbb{H}}$ (cf. the proof of
Lemma~\ref{lemma:Vanya-Kostya-n-4} and references therein). We
choose $\beta$ so that $W$ corresponds to the standard
five-dimensional representation of $\SS_6$ twisted by the sign
representation. Another way to describe the choice of $\beta$ is
through introducing the action of $N\slash\tilde{\mathbb{H}}$ on
the space $W^{\prime}=\Lambda^2(\C^4)$ (see~\cite{Nie92}).}
$$
\xymatrix{&1\ar@{->}[rr]&&\tilde{\mathbb{H}}\ar@{->}[rr]^{\alpha}&&N\ar@{->}[rr]^{\beta}&&\SS_{6}\ar@{->}[rr]&&1,}%
$$
where
$\tilde{\mathbb{H}}=\langle\mathbb{H},\mathrm{diag}(\sqrt{-1})\rangle$.
One can show that $N$ is a primitive subgroup of
$\SL_4(\mathbb{C})$.
\end{example}

The following theorem provides an explicit list of possible finite
subgroups in $\PGL_4(\mathbb{C})$ that satisfy all
hypothesis of Theorem~\ref{theorem:Vanya-Kostya}

\begin{theorem}[{see \cite[Chapter~VII]{Bli17} or \cite[\S
8.5]{Fe71}}] \label{theorem:Blichfeldt} Let $G$ be a primitive
subgroup of $\SL_4(\mathbb{C})$ such that
$Z(G)\subseteq[G,G]$. Then one of the following possibilities
holds:
\begin{itemize}
\item either $G$ satisfies the~hypotheses of Lemma~\ref{lemma:Segre} or Lemma~\ref{lemma:Segre-2},%

\item or $G$ is one of the~following groups:
\begin{itemize}
\item $\A_5$ or $\SS_5$,
\item $\SL_2(\mathbb{F}_{5})$,%
\item $\SL_2(\mathbb{F}_{7})$,%
\item $2.\A_{6}$, which is a~central extension of the~group $\A_{6}\cong\bar{G}$,%
\item $2.\SS_{6}$, which is a~central
extension\footnote{There are three non-isomorphic non-trivial
central extensions of the group $\SS_{6}$ with the~center isomorphic
to~$\mathbb{Z}_{2}$, two of which are embedded in
$\SL_4(\mathbb{C})$ (cf. \cite{Atlas}). But up to
conjugation there is only one subgroup of $\PGL_4(\mathbb{C})$
isomorphic to $\SS_6$.}
of the~group $\SS_{6}\cong\bar{G}$,%

\item $2.\A_{7}$, which is a~central extension of the~group $\A_{7}\cong\bar{G}$,%
\item $\Sp$,%
\item in the~notation of  Example~\ref{example:dim-4-primitive},
a~primitive subgroup in $N$ that contains~\mbox{$\alpha(\tilde{\mathbb{H}})$}.%
\end{itemize}
\end{itemize}
\end{theorem}

It should be pointed out that Theorem~\ref{theorem:Blichfeldt}
describes primitive subgroups of $\SL_4(\mathbb{C})$ up to
conjugation. Namely, if there are two monomorphisms
\mbox{$\iota_{1}\colon G\to\SL_4(\mathbb{C})$} and
\mbox{$\iota_{2}\colon G\to\SL_4(\mathbb{C})$} such that both
subgroups $\iota_{1}(G)$ and $\iota_{2}(G)$ are primitive, then it
follows from \cite[Chapter~VII]{Bli17} that $\iota_{1}(G)$ and
$\iota_{2}(G)$ are conjugate, but it may happen that
the representations of the group~$G$ given by $\iota_1$ and $\iota_2$ are
non-isomorphic, i.\,e.
there is no
element $g\in\SL_4(\mathbb{C})$ that makes the~diagram
$$
\xymatrix{
&G\ar@{->}[dl]_{\iota_1}\ar@{->}[dr]^{\iota_2}&\\
\SL_4\Big(\mathbb{C}\Big)\ar@{->}[rr]^{\theta_g}&&\SL_4\Big(\mathbb{C}\Big)\\
}%
$$
commutative, where $\theta_g$ is the~conjugation by $g$
(cf.~\cite{Atlas}).

\begin{lemma}
\label{lemma:A6} Suppose that $G\cong2.\A_6$. Then $G$ has no
semi-invariants of degree at most $4$.
\end{lemma}

\begin{proof}
Semi-invariants of $G$ are its~invariants by
Remark~\ref{remark:semiinvariants-vs-invariants}, and $G$ has no
odd degree in\-va\-ri\-ants, because $G$ contains a~scalar matrix
whose non-zero entries~are~$-1$.

To complete the~proof, it is enough to prove that $G$ has no invariants
of degree~$4$.

Let $V\cong\mathbb{C}^{4}$ be the~irreducible representation of
the~group $G$ that corresponds to the~embedding
$G\subset\SL_4(\mathbb{C})$. Without loss of generality,
we may assume that~\mbox{$\Lambda^{2}V\cong\mathbb{C}^{6}$} is
a~permutation representation of the~group $G\slash
Z(G)\cong\A_{6}$, because~$G$
has two four-dimen\-si\-o\-nal irreducible representations, which
give one subgroup $G\subset\SL_4(\mathbb{C})$ up to
conjugation.

Let $\chi$ be the character of the representation $V$,
and let $\chi_{4}$ be the~character
of the~representation~\mbox{$\mathrm{Sym}^{4}(V)$}. Then
$$
\chi_{4}\big(g\big)=\frac{1}{24}\Big(\chi\big(g\big)^4+6\chi\big(g\big)^2\chi\big(g^2\big)+3\chi\big(g^2\big)^2+8\chi\big(g\big)\chi\big(g^3\big)+6\chi\big(g^4\big)\Big)%
$$
for every $g\in G$. The~values of the~characters $\chi$ and $\chi_{4}$ are
listed in the~following table.
\begin{center}\renewcommand\arraystretch{1.3}
\begin{tabular}{|c|c|c|c|c|c|c|c|c|c|c|}
\hline
& $[5,1]_{10}$ & $[5,1]_5$ & $[4,2]_8$ & $[3,3]_6$ & $[3,3]_3$ &
$[3,1,1,1]_6$ & $[3,1,1,1]_3$ & $[2,2,1,1]_4$ & $z$ & $e$\\
\hline
\# & $144$ & $144$ & $180$ & $40$ & $40$ & $40$ & $40$ & $90$ & $1$ & $1$\\
\hline
$\chi$ & $1$ & $-1$ & $0$ & $-1$ & $1$ & $2$ & $-2$ & $0$ & $-4$ & $4$\\
\hline
$\chi_{4}$ & $0$ & $0$ & $-1$ & $2$ & $2$ & $-4$ & $-4$ & $3$ & $35$ & $35$\\
\hline
\end{tabular}
\end{center}
where the~first row lists the~types of the~elements in $G$ (for
example, the~symbol $[5,1]_{10}$ denotes the~set\footnote{ Note
that these sets do not coincide with conjugacy classes. For
example, the~image of the~set of the~ elements of type $[5,
1]_{10}$ under the~natural projection
$2.\A_6\to\A_6$ is a~union of two different
conjugacy classes in $\A_6$.} of order $10$ elements whose
image in $\A_6$ is a~product of disjoint cycles of
length~$5$~and~$1$), and $z$ and $e$ are the~non-trivial element
in the~center of $G$ and the~identity element, respectively.

Now one can check that the~inner product of the~character
$\chi_{4}$ and the~trivial character~is~zero, which implies that
the~subgroup $G$ does not have invariants of degree $4$.
\end{proof}

\begin{lemma}
\label{lemma:S6-A7} If $G\cong 2.\SS_6$ or $G\cong 2.\A_7$,
then $G$ has no semi-invariants of degree at most $4$.
\end{lemma}

\begin{proof}
Recall that these groups contain $2.\A_6$ and we can apply
Lemma~\ref{lemma:A6}.
\end{proof}

\begin{lemma}
\label{lemma:Vanya-Kostya-n-4}
Under the assumptions of Theorem~\ref{theorem:Blichfeldt}
the~subgroup $G$ has no
semi-invariants of degree at most $4$ if and only if $G$ is one of
the~following~groups:
\begin{itemize}
\item $2.\A_{6}$, $2.\SS_{6}$ or $2.\A_{7}$,%
\item $\Sp$,%
\item in the notation of Example~\ref{example:dim-4-primitive}, a
subgroup of $N$ that satisfies one of the following four
conditions:
\begin{itemize}
\item $G=N$,%

\item $\alpha(\tilde{\mathbb{H}})\subsetneq G$ and $\beta(G)\cong\A_{6}$,%

\item $\alpha(\tilde{\mathbb{H}})\subsetneq G$ and $\beta(G)\cong
\SS_{5}$, where the~embedding
$\beta(G)\subset\SS_{6}$ is non-standard, i.\,e.
the~standard one twisted by an~outer automorphism of~$\SS_6$,%

\item $\alpha(\tilde{\mathbb{H}})\subsetneq G$ and
$\beta(G)\cong\A_{5}$, where the~embedding
$\beta(G)\subset\SS_{6}$ is non-standard.%
\end{itemize}
\end{itemize}
\end{lemma}

\begin{proof}
Let $d$ be the~smallest positive number such $G$ has an
semi-invariant of degree~$d$. If $G\cong 2.\A_6$, then
$d\geqslant 5$ by Lemma~\ref{lemma:A6}. If $G\cong
2.\SS_{6}$ or $G\cong 2.\A_{7}$, then $d\geqslant
5$ by Lemma~\ref{lemma:S6-A7}. In fact, one can check by
direct computation that $d=8$ if $G\cong 2.\A_6$ or
$G\cong 2.\SS_{6}$ or $G\cong 2.\A_{7}$. If
$G\cong\SL_2(\mathbb{F}_{7})$, then the~equality $d=4$
holds by~\cite{MaSl73} and
Remark~\ref{remark:semiinvariants-vs-invariants}. If $G\cong\Sp$,
then the~equality $d=12$ holds by~\cite{Ma89}  and
Remark~\ref{remark:semiinvariants-vs-invariants}.

Suppose that $G\cong\SL_2(\mathbb{F}_{5})\cong
2.\A_{5}$. Then there is a~$\bar{G}$-invariant smooth
rational cubic curve $C\subset\mathbb{P}^{3}$, because
the~representation $G\to\GL_4(\mathbb{C})$ is a~symmetric
square of a~two-dimensional representation of the~group $G$.
The surface swept out by the lines tangent to the curve~$C$
is a $\bar{G}$-invariant surface of degree~$4$ (cf. proof
of Theorem~\ref{theorem:Vanya-Kostya-SL-4-weakly-exceptional}). Therefore, the
inequality $d\leqslant 4$ holds\footnote{Actually, one can show that $d=4$
in this case.}.

Let us use the~notation of Example~\ref{example:dim-4-primitive}.
By Theorem~\ref{theorem:Blichfeldt},
Remark~\ref{remark:real-representation} and
Lemmas~\ref{lemma:Segre}~and~\ref{lemma:Segre-2}, to complete
the~proof we may assume that $G$ is a~primitive subgroup in $N$
that contains $\alpha(\tilde{\mathbb{H}})$.

One can show that the~group $\tilde{\mathbb{H}}$ has no invariants
of degree less than $4$  and its invariants of degree $4$ form
a~five-dimensional vector space $W$ (see e.\,g.
\cite[Lemma~3.18]{Pr07}).

The group $\beta(G)$ naturally acts on $W$. Moreover, the~subgroup
$G$ has an invariant of degree $4$ if and only if
the~representation $W$ has a~one-dimensional subrepresentation of
the~group $\beta(G)$. On the other hand, it follows from
\cite{Nie92} that if $G=N$, then $W$ is an irreducible
representation of $\beta(G)=\SS_6$.

It follows from \cite[\S 123]{Bli17} that, up to conjugation, there exist
exactly $9$ possibilities for the~subgroup $G\subset N$ such that
$G$~is~primi\-tive. These possibilities are listed in the~following table:
\begin{center}\renewcommand\arraystretch{1.3}
\begin{center}
\begin{tabular}{|c|c|c|c|}
\hline
Label of the~group $G$ & $\beta(G)$ & Generators of the~subgroup $\beta(G)\subseteq\SS_{6}$ & Splitting type\\
\hline $13^{\circ}$ & $\mathbb{Z}_{5}$ &  $(24635)$ & $1,1,1,1,1$\\
\hline $14^{\circ}$ & $\mathbb{Z}_{5}\rtimes\mathbb{Z}_{2}$ & $(24635), (36)(45)$ & $1,2,2$ \\
\hline $15^{\circ}$ & $\mathbb{Z}_{5}\rtimes\mathbb{Z}_{4}$ & $(24635), (3465)$ & $1,2,2$ \\
\hline
$16^{\circ}$ & $\A_5$ &  $(24635), (34)(56)$ & $1,4$\\
\hline
$17^{\circ}$ & $\A_5$ &  $(24635), (12)(36)$ & $5$\\
\hline
$18^{\circ}$ & $\SS_5$ &  $(24635), (56)$ & $1,4$\\
\hline
$19^{\circ}$ & $\SS_5$ &  $(24635), (12)(34)(56)$ & $5$\\
\hline
$20^{\circ}$ & $\A_6$ & $(24635), (12)(34)$ & $5$\\
\hline
$21^{\circ}$ & $\SS_6$ &  $(24635), (12)$ & $5$\\
\hline
\end{tabular}
\end{center}
\end{center}
where the~first column lists the~labels of the~subgroup $G$
according to~\cite[\S 123]{Bli17} and the~last column lists
the~dimensions of the~irreducible $\beta(G)$-subrepresentations of~$W$.

Note that $\mathbb{H}\subset\tilde{\mathbb{H}}$ has no
semi-invariants of degree $3$, because $\mathbb{H}$ has no
invariants of degree $3$, the~center of the~group $\mathbb{H}$
coincides with its commutator and acts non-trivially on cubic
forms.

The subgroups of $N$ described in
Lemma~\ref{lemma:Vanya-Kostya-n-4} are the~subgroups $21^{\circ}$,
$20^{\circ}$, $19^{\circ}$, $17^{\circ}$, respectively. We see
that $d\leqslant 4$  if $G$ is the~subgroup $13^{\circ}$,
$14^{\circ}$, $15^{\circ}$, $16^{\circ}$~or~$18^{\circ}$. On the
other hand, if $G$ is the~subgroup $17^{\circ}$, $19^{\circ}$,
$20^{\circ}$ or $21^{\circ}$, then the~subgroup $G$ has neither
semi-invariants of degree less than $4$, nor invariants of degree
$4$. Let us prove that the~subgroup $17^{\circ}$ does not have
semi-invariants of degree $4$. Since the absence of
semi-invariants of degree $4$ implies the absence of
semi-invariants of degree $2$, this would imply that in the~case
when $G$ is the~subgroup $17^{\circ}$, $19^{\circ}$,
$20^{\circ}$~or~$21^{\circ}$ of the group~$N$ the
inequality~\mbox{$d\ge 5$} holds\footnote{In fact, one can check
by direct computation that $d=8$ if $G$ is the~subgroup
$17^{\circ}$, $19^{\circ}$, $20^{\circ}$ or $21^{\circ}$.}.

Suppose that $G$ is the~subgroup $17^{\circ}$, and suppose, in
addition, that $G$ does have a~semi-invariant $\Phi$ of degree
$4$. Let us show that this assumption leads to a~contradiction.

Note that the polynomial $\Phi$ is not $\tilde{\mathbb{H}}$-invariant, because
$\Phi$ is not $G$-invariant and
\mbox{$G\slash\tilde{\mathbb{H}}\cong\beta(G)\cong\A_{5}$} is
a~simple group. Let $Z$ be the~center of the~group
$\tilde{\mathbb{H}}$. Put
$\bar{\mathbb{H}}=\phi(\tilde{\mathbb{H}})$. Then
\mbox{$\tilde{\mathbb{H}}\slash
Z\cong\bar{\mathbb{H}}\cong\mathbb{Z}_2^4$}, and $Z$ acts trivially
on $\Phi$. Thus, there is a~homomorphism
$\xi\colon\bar{\mathbb{H}}\to\mathbb{C}^{*}$ such that
$\mathrm{ker}(\xi)\ne\bar{\mathbb{H}}$, which implies that
$\mathrm{ker}(\xi)\cong\mathbb{Z}_{2}^{3}$, because
$\mathrm{im}(\chi)$ is a~cyclic group. Let
$\theta\colon\bar{G}\to\mathrm{Aut}(\bar{\mathbb{H}})$ be
the~homomorphism such that
$$
\theta\big(g\big)\Big(h\Big)=ghg^{-1}\in\bar{\mathbb{H}}\cong\mathbb{Z}_{2}^{4}
$$
for all $g\in\bar{G}$ and $h\in\bar{\mathbb{H}}$. Consider
$\bar{\mathbb{H}}$ as a~vector space over $\mathbb{F}_2$. Then
$\theta$ induces a~monomorphism
$\tau\colon\beta(G)\to\GL_4(\mathbb{F}_{2})$ and
$\mathrm{ker}(\xi)$ is a~$\mathrm{im}(\tau)$-invariant subspace.
But $\mathrm{im}(\tau)\cong\A_{5}$ has no non-trivial
three-dimensional representations over $\mathbb{F}_{2}$, because
$|\GL_3(\mathbb{F}_{2})|=168$ is not divisible by
$|\A_{5}|=60$. Thus, we see that there is a~non-zero
element $t\in\bar{\mathbb{H}}$ such that $t$ is
$\mathrm{im}(\tau)$-invariant. Let $F$ be the~stabilizer of $t$ in
$\GL_4(\mathbb{F}_{2})$. Then
$\A_{5}\cong\mathrm{im}(\tau)\subset F$, which is
impossible, because $|F|=1344$ is not divisible by
$|\A_{5}|=60$.
\end{proof}

\medskip

Combining the previous results we obtain the following.

\begin{theorem}
\label{theorem:Vanya-Kostya-n-4} Let $G$ be a finite subgroup in
$\SL_4(\C)$. Then the~following conditions are equivalent:
\begin{itemize}
\item the~singularity $(V\ni O)$ is exceptional,%

\item the~inequality $\mathrm{lct}(\mathbb{P}^{3},\bar{G})\geqslant 5/4$ holds,%

\item the~group $G$ is primitive and $G$ does not have semi-invariants of degree at most $4$,%

\item $\bar{G}=\phi(G')$, where $G'$ is one of the $8$
subgroups listed in Lemma~\ref{lemma:Vanya-Kostya-n-4}.
\end{itemize}
\end{theorem}

\begin{proof}
The required assertion follows from
Theorems~\ref{theorem:criterion}, \ref{theorem:Vanya-Kostya} and
\ref{theorem:Blichfeldt} and Lemma~\ref{lemma:Vanya-Kostya-n-4}.
\end{proof}

\section{Five-dimensional case}
\label{section:5}

The purpose of this section is to present an explicit
classification of exceptional five-dimensional singularities (see
Theorem~\ref{theorem:Vanya-Kostya-n-5},
cf.~Theorems~\ref{theorem:Shokurov-n-2},
\ref{theorem:Dima-Yura-n-3} and
Theorem~\ref{theorem:Vanya-Kostya-n-4}), and prove some relevant
results.

Let $\bar{G}$ be a~finite subgroup in
$\mathrm{Aut}(\mathbb{P}^{4})$, and consider the natural
projection
\mbox{$\phi\colon\SL_5(\mathbb{C})\to\mathrm{Aut}(\mathbb{P}^{4})\cong
\PGL_5(\mathbb{C})$}. Then there is a~finite subgroup
$G\subset\SL_5(\mathbb{C})$~such~that $\phi(G)=\bar{G}$.
Suppose that $G$ is primitive. Then we may assume that
$Z(G)\subseteq[G,G]$ (see \cite{Br67} and \cite{Fe71}).

\begin{example}[{cf. Appendix~\ref{section:HM}}]
\label{example:dim-5-primitive} Let $\mathbb{H}$ be the~Heisenberg
group of all unipotent $3\times 3$-matrices
with~entries~in~$\mathbb{F}_{5}$. Then there is a~monomorphism
$\rho\colon\mathbb{H}\to\SL_5(\mathbb{C})$. Let
$\mathbb{HM}$ be the~normalizer of the~subgroup $\rho(\mathbb{H})$
$\subset\SL_5(\mathbb{C})$. Then there is an~exact
sequence
$$
\xymatrix{&1\ar@{->}[rr]&&\mathbb{H}\ar@{->}[rr]^{\alpha}&&\mathbb{HM}\ar@{->}[rr]^{\beta}&&\SL_2\Big(\mathbb{F}_{5}\Big)\ar@{->}[rr]&&1,}%
$$
and $\mathbb{HM}$ is a primitive subgroup in
$\SL_5(\mathbb{C})$ (see \cite[Theorem~9A]{Br67},
\cite{HoMu73}).
\end{example}

\begin{theorem}[{see \cite{Br67} or \cite[\S 8.5]{Fe71}}]
\label{theorem:Brauer} Let $G$ be a finite primitive subgroup in
$\SL_5(\mathbb{C})$ such that $Z(G)\subseteq[G,G]$. Then
$G$ is one of the groups $\A_5$,
$\A_6$, $\SS_5$, $\SS_6$,
$\PSL_2(\mathbb{F}_{11})$, $\PSp$, or, in the~notation of
Example~\ref{example:dim-5-primitive}, a~primitive subgroup of
$\mathbb{HM}$ that contains $\alpha(\mathbb{H})$.
\end{theorem}

Note that if there are two monomorphisms $\iota_{1}\colon
G\to\SL_5(\mathbb{C})$ and \mbox{$\iota_{2}\colon
G\to\SL_5(\mathbb{C})$} such that both subgroups
$\iota_{1}(G)$ and $\iota_{2}(G)$ are primitive, then
$\iota_{1}(G)$ and $\iota_{2}(G)$ are conjugate.

\begin{lemma}
\label{lemma:n-5-1} Suppose that $G$ is one of the~following
groups: $\A_5$, $\A_6$, $\SS_5$,
$\SS_6$, $\PSL_2(\mathbb{F}_{11})$ or $\PSp$. Then
$G$ has an invariant of degree at most $4$, which implies that
$\mathrm{lct}(\mathbb{P}^{4},\bar{G})\leqslant 4/5$.
\end{lemma}

\begin{proof}
If $G$ is $\A_5$, $\A_6$,
$\SS_5$ or~$\SS_6$, then $G$ has an invariant of
degree $2$ by Remark~\ref{remark:real-representation}. If
$G\cong\PSp$, then $G$ has an invariant of degree $4$ (see
\cite{Bu91}).  If $G\cong\PSL_2(\mathbb{F}_{11})$, then
$G$ has an invariant of degree $3$ (see \cite{Ad78}).
\end{proof}

\begin{lemma}
\label{lemma:lemma:n-5-4} In the~notation of
Example~\ref{example:dim-5-primitive}, suppose that
$\alpha(\mathbb{H})\subsetneq G\subseteq\mathbb{HM}$. Then~$G$
has no semi-invariants of degree at most $5$ if and only if
 either $G=\mathbb{HM}$ or~$G$ is a~subgroup of $\mathbb{HM}$ of index~ $5$.
\end{lemma}

\begin{proof}
Let $V$ be the~vector space of $\mathbb{H}$-invariant forms of
degree $5$. Then the~group
$\mathbb{HM}\slash\alpha(\mathbb{H})\cong\SL_2(\mathbb{F}_{5})\cong
2.\A_5$ naturally acts on the~vector space $V$. Moreover,
it follows from \cite[Theorem~3.5]{HoMu73} that
$V=V^{\prime}\oplus V^{\prime\prime}$, where $V^{\prime}$ and
$V^{\prime\prime}$ are three-dimensional
$\mathrm{im}(\beta)$-invariant linear subspaces that arise from
two~non-equi\-va\-lent three-dimensional representations of
the~group $\A_5$, respectively. Therefore, we see that $G$
has a~semi-invariant of degree $5$ if and only if $V^{\prime}$ has
a~$\beta(G)$-invariant one-dimensional subspace.

Let $Z\cong\mathbb{Z}_2$ be the~center of the~group
$\mathbb{HM}\slash\alpha(\mathbb{H})\cong 2.\A_5$. Then
$2.\A_5/Z\cong\A_{5}$. Moreover, either $\beta(G)$
is cyclic, or $Z\subseteq \beta(G)$ and $\beta(G)/Z$ is one of
the~following subgroups of $\A_5$: dihedral group of order
$6$, dihedral group of order $10$, the~group
$\mathbb{Z}_{2}\times\mathbb{Z}_{2}$, the~group $\A_4$,
the~group $\A_5$.

If $\beta(G)$ is cyclic, then $V^{\prime}$ is a~sum of
one-dimensional $\beta(G)$-invariant linear subspaces. Hence we
may assume that $Z\subseteq \beta(G)$. Recall that
$Z\cong\mathbb{Z}_{2}$ acts trivially on~$V^{\prime}$. Thus,  if
$\beta(G)\slash Z\cong\mathbb{Z}_{2}\times\mathbb{Z}_{2}$, then
$V^{\prime}$ is a~sum of one-dimensional $\beta(G)$-invariant
subspaces.

If $\beta(G)/Z$ is a~dihedral group, then $V^{\prime}$ must have
one-dimensional $\beta(G)$-invariant subspace, because irreducible
representations of dihedral groups are one-dimensional or two-dimensional.

If $\beta(G)/Z\cong\A_5$ or $\beta(G)/Z\cong\A_4$,
then $V^{\prime}$ is an irreducible representation
of~\mbox{$\beta(G)/Z$}, which implies that  $V^{\prime}$ is an irreducible
representation of the~group $\beta(G)$. Now using
Corollary~\ref{corollary:Heisenberg-representations}, we complete
the~proof.
\end{proof}

The main purpose of this section is to prove the~following result.

\begin{theorem}
\label{theorem:HM} In the~notation of
Example~\ref{example:dim-5-primitive}, let $G$ be a~subgroup of
the~group~$\mathbb{HM}$ of index $5$. Then
$\mathrm{lct}(\mathbb{P}^{4},\bar{G})\geqslant 6/5$.
\end{theorem}

Combining the previous results we obtain the following.

\begin{theorem}
\label{theorem:Vanya-Kostya-n-5} Let $G$ be a finite subgroup in
$\SL_5(\C)$. Then the~following conditions are~equivalent:
\begin{itemize}
\item the~singularity $(V\ni O)$ is exceptional,%

\item the~inequality $\mathrm{lct}(\mathbb{P}^{4},\bar{G})\geqslant 6/5$ holds,%

\item the~group $G$ is primitive
and $G$ does not have semi-invariants of degree at most $5$,%

\item in the~notation of Example~\ref{example:dim-5-primitive},
either $G\cong\mathbb{HM}$ or $G$ is isomorphic to a subgroup of
the group $\mathbb{HM}$ of index $5$.
\end{itemize}
\end{theorem}

\begin{proof}
The required assertion follows from
Theorems~\ref{theorem:criterion}, \ref{theorem:HM},
\ref{theorem:Brauer} and Lemmas~\ref{lemma:lemma:n-5-4} and
\ref{lemma:n-5-1}.
\end{proof}

In the remaining part of this section we will prove
Theorem~\ref{theorem:HM}. Let us use the notation of
Example~\ref{example:dim-5-primitive}. Suppose that $G$ be
a~subgroup of the~group $\mathbb{HM}$ of index~$5$.

\begin{lemma}
\label{lemma:A4-in-HM-orbits}  Let $\Lambda$ be
a~$\bar{G}$-invariant subset of $\mathbb{P}^4$. Then $\Lambda$
consists of at least $10$ points.
\end{lemma}

\begin{proof}
The required assertion follows from Lemma~\ref{lemma:lemma:n-5-4}
and Corollary~\ref{corollary:Heisenberg-representations}.
\end{proof}

Suppose that~$\mathrm{lct}(\mathbb{P}^{4},\bar{G})<6/5$. Let us
derive a contradiction.

There is a~rational positive number $\lambda<6/5$ and an~effective
$\bar{G}$-invariant \mbox{$\mathbb{Q}$-divisor}~$D$ on
$\mathbb{P}^{5}$~such~that $D\sim_{\mathbb{Q}}-K_{\mathbb{P}^{4}}$
and the~log pair $(\mathbb{P}^{4},\lambda D)$ is strictly log
canonical. Let~$S$ be an irreducible subvariety of
$\mathbb{P}^{4}$ that is a~minimal center in
$\mathbb{LCS}(\mathbb{P}^{4},\lambda D)$. Then $S$ is either a
curve or a surface by
Lemma~\ref{lemma:semiinvariants-and-codimension-one}.

Let $Z$ be the~$\bar{G}$-orbit of the~subvariety
$S\subset\mathbb{P}^{4}$, and let $r$ be the~number of irreducible
components of the~subvariety $Z$. We may assume that
$$
\mathbb{LCS}\Big(\mathbb{P}^{4},\lambda D\Big)=\bigcup_{g\in\bar{G}}\Big\{g\big(S\big)\Big\}%
$$
by Lemma~\ref{lemma:Kawamata-Shokurov-trick}. Then
$\mathrm{Supp}(Z)=\mathrm{LCS}(\mathbb{P}^{4},\lambda D)$. It
follows from Lemma~\ref{lemma:centers} that
$$
S\cap g\big(S\big)\ne \varnothing\iff S=g\big(S\big)%
$$
for every $g\in\bar{G}$. Then $\mathrm{deg}(Z)=r\mathrm{deg}(S)$.

Let $\mathcal{I}$ be the~multiplier ideal sheaf of the~log pair
$(\mathbb{P}^{4},\lambda D)$, and let $\mathcal{L}$ be the~log
canonical singularities subscheme of the~log pair
$(\mathbb{P}^{4},\lambda D)$. By
Theorem~\ref{theorem:Shokurov-vanishing}, there is an~exact
sequence
\begin{equation}
\label{equation:exact-sequence-n-5}
0\longrightarrow H^{0}\Big(\mathcal{O}_{\mathbb{P}^{4}}\big(n\big)\otimes\mathcal{I}\Big)\longrightarrow H^{0}\Big(\mathcal{O}_{\mathbb{P}^{4}}\big(n\big)\Big)\longrightarrow H^{0}\Big(\mathcal{O}_{\mathcal{L}}\otimes\mathcal{O}_{\mathbb{P}^{4}}\big(n\big)\Big)\longrightarrow 0%
\end{equation}
for every $n\geqslant 1$. Note that $Z=\mathcal{L}$
by Remark~\ref{remark:log-canonical-subscheme}.

\begin{lemma}
\label{lemma:n-5-curves} The center $S$ is not a curve.
\end{lemma}

\begin{proof}
Suppose that $S$ is a curve.  Then it follows from
Theorem~\ref{theorem:Kawamata} that $S$ is a~smooth curve of genus
$g$ such that  $2g-2<\mathrm{deg}(S)$. Moreover, it follows from
Lemma~\ref{lemma:degree} that $\mathrm{deg}(Z)\leqslant 10$.~Then
$2g-2<\mathrm{deg}(S)\leqslant 10$, which implies that~\mbox{$g\leqslant 5$}.
The curve $Z$ is not contained in a~hyperplane, because $G$ is
transitive. Then \mbox{$10\geqslant \mathrm{deg}(Z)=r\mathrm{deg}(S)$},
which implies that $r\leqslant 10$.

Using $(\ref{equation:exact-sequence-n-5})$ and the~Riemann--Roch
theorem, we see~that
\begin{equation}
\label{equation:inequalities-n-5-curves} 5=h^{0}\Big(\mathcal{O}_{\mathcal{L}}\otimes\mathcal{O}_{\mathbb{P}^{3}}\big(1\big)\Big)=r\Big(\mathrm{deg}\big(S\big)-g+1\Big),%
\end{equation}
because $\mathcal{L}=Z$ and $2g-2<\mathrm{deg}(S)$. Thus, either
$r=1$ or $r=5$.

If $r=5$, then $\mathrm{deg}(S)=2$ and $g=0$, which contradicts
$(\ref{equation:inequalities-n-5-curves})$. We see that $r=1$.
Thus $S$ is a $\bar{G}$-invariant irreducible curve of genus $g\le 5$, which
is impossible by Lemma~\ref{lemma:A4-in-HM-elliptic}.
\end{proof}

We see that $S$ is a surface. Then $\mathrm{deg}(Z)\leqslant 10$
by Lemma~\ref{lemma:degree}. It follows from
Theorem~\ref{theorem:Kawamata}~that $S$ is normal and has at most
rational singularities, and there is an~effective
$\mathbb{Q}$-divisor $B_{S}$ and an~ample $\mathbb{Q}$-divisor
$\Delta$ on the~surface $S$ such that
$$
K_{S}+B_{S}+\Delta\qlin\mathcal{O}_{\mathbb{P}^{4}}\big(1\big)\Big\vert_{S},%
$$
and the~log pair $(S,B_{S})$ has Kawamata log terminal
singularities. Therefore,
the equality $r=1$ holds, since
two irreducible surfaces in $\mathbb{P}^{4}$ have non-empty
intersection.

Thus, we see that the~surface $S=Z$ is $\bar{G}$-invariant.

\begin{lemma}
\label{lemma:n-5-surface-plane} The surface $S$ is not contained
in a~hyperplane in $\mathbb{P}^{4}$.
\end{lemma}

\begin{proof}
The required assertion follows from the~fact that $G$ is transitive.
\end{proof}

\begin{lemma}
\label{lemma:n-5-surface-quadric} The surface $S$ is not contained
in a~quadric hypersurface in $\mathbb{P}^{4}$.
\end{lemma}

\begin{proof}
Suppose that there is a~quadric hypersurface
$Q\subset\mathbb{P}^{4}$ such that $S\subset Q$. Then~$Q$ is
irreducible by Lemma~\ref{lemma:n-5-surface-plane}. Moreover, it
follows from Lemma~\ref{lemma:lemma:n-5-4} that there is a~quadric
hypersurface $Q^{\prime}\subset\mathbb{P}^{4}$ such that
$S\subseteq Q\cap Q^{\prime}$, because otherwise the~quadric~$Q$
would be $\bar{G}$-invariant. Then $Q^{\prime}$ is irreducible by
Lemma~\ref{lemma:n-5-surface-plane}.

Suppose that $S=Q\cap Q^{\prime}$.
If $S$ is non-singular, consider a~pencil $\mathcal{P}$
generated by the quadrics $Q$ and~$Q'$.
Then $\mathcal{P}$ contains exactly $5$ singular quadrics, which
are simple quadric cones. This means that there is a
$\bar{G}$-orbit in $\mathbb{P}^{4}$ consisting of at~most~$5$
points, which is impossible, because $G$ has no semi-invariants of
degree up to~$5$.
Therefore, the surface $S$ is singular.

It follows
from~\cite{HiWa81}~that $|\mathrm{Sing}(S)|\leqslant 4$, because
$S$ has canonical singularities since~$S$ is a~complete
intersection that has Kawamata log terminal singularities. But
$\mathrm{Sing}(S)$ is $\bar{G}$-invariant, which contradicts
Lemma~\ref{lemma:A4-in-HM-orbits}.

We see that $S\ne Q\cap Q^{\prime}$. Therefore, it follows from
Lemma~\ref{lemma:n-5-surface-plane} that either $S$ is a~cone over
a~smooth rational cubic curve, or $S$ is a~smooth cubic scroll.

If $S$ is a~cone, then its vertex is $\bar{G}$-invariant, which
is impossible since $G$ is transitive. Thus, we see that
$S$ is a~smooth cubic scroll. Then there is a~unique line
$L\subset S$ such that $L^{2}=-1$, which implies that $L$ must be
$\bar{G}$-invariant, which is again impossible, because $G$ is
transitive.
\end{proof}

Let $H$ be a~hyperplane section of the~surface
$S\subset\mathbb{P}^{4}$.

\begin{lemma}
\label{lemma:H-invariant-irreducible-surface} The equalities
$H\cdot H=-H\cdot K_{S}=5$ and $\chi(\mathcal{O}_{S})=0$ hold.
\end{lemma}

\begin{proof}
It follows from
Corollary~\ref{corollary:Heisenberg-representations} that there is
$m\geqslant 0$ such that one has
\mbox{$h^{0}(\mathcal{O}_{\mathbb{P}^{4}}(3)\otimes\mathcal{I})=5m$}. Let
us show that this is possible only if $H\cdot H=-H\cdot K_{S}=5$
and $\chi(\mathcal{O}_{S})=0$.

It follows from the~Riemann--Roch theorem and
Theorem~\ref{theorem:Shokurov-vanishing} that
\begin{equation}
\label{equation:RR}
h^{0}\Big(\mathcal{O}_{S}\big(nH\big)\Big)=\chi\Big(\mathcal{O}_{S}\big(nH\big)\Big)=\chi\big(\mathcal{O}_{S}\big)+\frac{n^{2}}{2}\Big(H\cdot H\Big)-\frac{n}{2}\Big(H\cdot K_{S}\Big)%
\end{equation}
for any $n\geqslant 1$. It follows from
Lemma~\ref{lemma:n-5-surface-plane},
the~equalities~\ref{equation:RR} and the~exact
sequence~\ref{equation:exact-sequence-n-5}~that
\begin{equation}
\label{equation:RR-1}
5=h^{0}\Big(\mathcal{O}_{S}\big(H\big)\Big)=\chi\big(\mathcal{O}_{S}\big)+\frac{1}{2}\Big(H\cdot H\Big)-\frac{1}{2}\Big(H\cdot K_{S}\Big),%
\end{equation}
and it follows from Lemma~\ref{lemma:n-5-surface-quadric},
the~equalities~\ref{equation:RR}  and the~exact
sequence~\ref{equation:exact-sequence-n-5} that
\begin{equation}
\label{equation:RR-2}
15=h^{0}\Big(\mathcal{O}_{S}\big(2H\big)\Big)=\chi\big(\mathcal{O}_{S}\big)+2\Big(H\cdot H\Big)-\Big(H\cdot K_{S}\Big).%
\end{equation}

It follows from
Lemmas~\ref{lemma:degree},~\ref{lemma:n-5-surface-plane} and
\ref{lemma:n-5-surface-quadric} that $4\leqslant H\cdot
H=\mathrm{deg}(S)\leqslant 10$.

Suppose that $H\cdot H=10$. It follows from
the~equalities~\ref{equation:RR-1} and \ref{equation:RR-2} that
$\chi(\mathcal{O}_{S})=5$~and~ $H\cdot K_{S}=H\cdot H=10$, which
is impossible, because $H\qlin K_{S}+B_{S}+\Delta$, where
$\Delta$ is ample and $B_{S}$ is effective.
Thus $H\cdot H\le 9$.

It follows from the~equalities~\ref{equation:RR-1} and
\ref{equation:RR-2} that
$$
H\cdot K_{S}=3\chi\big(\mathcal{O}_{S}\big)-5=3\Big(H\cdot H\Big)-20.%
$$

It follows from the~equalities~\ref{equation:RR} and the~exact
sequence~\ref{equation:exact-sequence-n-5} that
$$
h^{0}\Big(\mathcal{O}_{\mathbb{P}^{4}}\big(3\big)\otimes\mathcal{I}\Big)=35-h^{0}\Big(\mathcal{O}_{S}\big(3H\big)\Big)=35-\Bigg(\chi\big(\mathcal{O}_{S}\big)+\frac{9}{2}\Big(H\cdot H\Big)-\frac{3}{2}\Big(H\cdot K_{S}\Big)\Bigg)=5m,%
$$
which implies that $H\cdot H=5$, $\chi(\mathcal{O}_{S})=0$ and
$H\cdot K_{S}=-5$, because  $4\leqslant H\cdot H\leqslant 9$.
\end{proof}

Let $\pi\colon U\to S$ be the~minimal resolution of the~surface
$S$. Then $\kappa(U)=-\infty$ and
$$
1-h^{1}\big(\mathcal{O}_{U}\big)=1-h^{1}\big(\mathcal{O}_{S}\big)=h^{2}\big(\mathcal{O}_{S}\big)=h^{2}\big(\mathcal{O}_{U}\big)=h^{0}\Big(\mathcal{O}_{U}\big(K_{U}\big)\Big)=0,
$$
because $S$ has rational singularities and $\kappa(U)=-\infty$ since
$H\cdot K_{S}=-5<0$.

\begin{corollary}
\label{corollary:n-5-elliptic-surface} The surface $S$ is
birational to $E\times\mathbb{P}^{1}$, where $E$ is smooth
elliptic curve.
\end{corollary}

By Remark~\ref{remark:stupid-action}, there is a~monorphism
$\xi\colon\bar{G}\to\mathrm{Aut}(Y)$, which contradicts
Corollary~\ref{corollary:elliptic-surface}.

The obtained contradiction completes the~proof of
Theorem~\ref{theorem:HM}.

\appendix

\section{Horrocks--Mumford group}
\label{section:HM}

Let $\mathbb{H}$ be the~Heisenberg group of all unipotent $3\times
3$-matrices~with~entries~in~$\mathbb{F}_{5}$. Then
$$
\mathbb{H}=\Big\langle x, y, z\ \Big|\ x^5=y^5=z^5=1,\ xz=zx,\ yz=zy,\ xy=zyx\Big\rangle%
$$
for some $x,y,z\in\mathbb{H}$. There is
a~monomorphism
$\rho\colon\mathbb{H}\to\SL_5(\mathbb{C})$
such that
$$
\rho\big(x\big)=\left(
\begin{array}{ccccc}
0 & 0 & 0& 0 & 1\\
1 & 0 & 0& 0 & 0\\
0 & 1 & 0& 0 & 0\\
0 & 0 & 1& 0 & 0\\
0 & 0 & 0& 1 & 0\\
\end{array}
\right),\ \rho\big(y\big)=\left(
\begin{array}{ccccc}
\zeta & 0 & 0& 0 & 0\\
0 & \zeta^2 & 0& 0 & 0\\
0 & 0 & \zeta^3& 0 & 0\\
0 & 0 & 0& \zeta^4 & 0\\
0 & 0 & 0& 0 & 1\\
\end{array}
\right),
$$
where $\zeta$  is a~non-trivial fifth root of unity. Let us
identify $\mathbb{H}$ with $\mathrm{im}(\rho)$. Then
$Z(\mathbb{H})\cong \mathbb{Z}_{5}$ and
$$
\left(\begin{array}{ccccc}
\zeta & 0 & 0& 0 & 0\\
0 & \zeta & 0& 0 & 0\\
0 & 0 & \zeta& 0 & 0\\
0 & 0 & 0& \zeta & 0\\
0 & 0 & 0& 0 & \zeta\\
\end{array}
\right)\in Z\Big(\mathbb{H}\Big),
$$
where $Z(\mathbb{H})$ is the~center of $\mathbb{H}$. Let
$\phi\colon\GL_5(\mathbb{C})\to\PGL_5(\mathbb{C})$
be the~natural projection.

\begin{lemma}[{\cite[\S 1]{HoMu73}}]
\label{lemma:Heisenberg-representations} Let
$\chi\colon\mathbb{H}\to\GL_N(\mathbb{C})$ be an
irreducible representation of $\mathbb{H}$. Then either $N=1$ and
$Z(\mathbb{H})\subseteq\mathrm{ker}(\chi)$, or $N$ is divisible by
$5$.
\end{lemma}

Take $n\in\mathbb{Z}_{\geqslant 0}$.
Then $\mathbb{H}$ naturally acts on
$H^{0}(\mathcal{O}_{\mathbb{P}^{4}}(n))$.

\begin{corollary}
\label{corollary:Heisenberg-representations} Let $V$ be
a~$\mathbb{H}$-invariant subspace in
$H^{0}(\mathcal{O}_{\mathbb{P}^{4}}(n))$. Then either
$\mathrm{dim}(V)$ is divisible by $5$, or $n$ is divisible by $5$.
\end{corollary}

Let $\mathbb{HM}\subset\SL_5(\mathbb{C})$ be
the~normalizer of the~subgroup $\mathbb{H}$. Then there is
an~exact sequence
$$
\xymatrix{&1\ar@{->}[rr]&&\mathbb{H}\ar@{->}[rr]^{\alpha}&&\mathbb{HM}\ar@{->}[rr]^{\beta}&&\SL_2\Big(\mathbb{F}_{5}\Big)\ar@{->}[rr]&&1,}%
$$
and it follows from~\cite[\S 1]{HoMu73} that there is a~subgroup
$\mathbb{M}\subset\mathbb{HM}$ such that
$\mathbb{HM}=\mathbb{H}\rtimes\mathbb{M}$ and
$\mathbb{M}\cong\beta(\mathbb{M})=\SL_2(\mathbb{F}_{5})\cong
2.\A_{5}$. Put $\bar{\mathbb{H}}=\phi(\mathbb{H})$ and
$\overline{\mathbb{HM}}=\phi(\mathbb{HM})$. Then
\mbox{$\overline{\mathbb{HM}}\slash\bar{\mathbb{H}}\cong\SL_2(\mathbb{F}_{5})$}
and $\bar{\mathbb{H}}\cong\mathbb{Z}_{5}\times\mathbb{Z}_{5}$. Let
$Z(\mathbb{HM})$ be the~center of the~group $\mathbb{HM}$. Then
$Z(\mathbb{HM})=Z(\mathbb{H})\cong \mathbb{Z}_{5}$.

\begin{corollary}
\label{corollary:HM-center} The group $\overline{\mathbb{HM}}$ is
isomorphic to $\mathbb{HM}\slash Z(\mathbb{HM})$.
\end{corollary}

Let $G$ be a~subgroup of the~group $\mathbb{HM}$ of index $5$.
Then $G\cong\mathbb{H}\rtimes
2.\A_{4}\subset\mathbb{H}\rtimes 2.\A_{5}$ and
$|\bar{G}|=600$, where $\bar{G}=\phi(G)$. Let $Z(G)$ be the~center
of the~group $G$. Then
$Z(G)=Z(\mathbb{HM})=Z(\mathbb{H})\cong\mathbb{Z}_{5}$.

\begin{lemma}
\label{lemma:A4-non-trivial-action} Let $g$ be an~element of the~
group $\bar{G}$ such that $gh=hg\in\bar{G}$ for every element
$h\in\bar{\mathbb{H}}$. Then $g\in\bar{\mathbb{H}}$.
\end{lemma}
\begin{proof}
The required assertion follows from~\cite[\S 1]{HoMu73}.
\end{proof}

\begin{lemma}
\label{lemma:normal-subgroups-in-A4} Let $F$ be a~proper normal
subgroup of $2.\A_4$. Then either $F\cong\mathbb{Z}_{2}$
is a~center of the~group $2.\A_4$, or
$F\cong\mathbb{Q}_{8}$, where $\mathbb{Q}_{8}$ is the~quaternion
group of order $8$.
\end{lemma}

\begin{proof}
The only nontrivial normal subgroup of the~group $\A_4$ is
isomorphic to the~ group~\mbox{$\mathbb{Z}_{2}\times\mathbb{Z}_{2}$}.
\end{proof}

\begin{lemma}
\label{lemma:A4-active-action} The group $\bar{\mathbb{H}}$ contains no proper
non-trivial subgroups that are normal in $\bar{G}$.
\end{lemma}

\begin{proof}
Let $\theta\colon \overline{\mathbb{HM}}\to\mathrm{Aut}(\bar{\mathbb{H}})$ be
the~homomorphism such that
$$
\theta\big(g\big)\Big(h\Big)=ghg^{-1}\in\bar{\mathbb{H}}
$$
for all $g\in \overline{\mathbb{HM}}$ and $h\in\bar{\mathbb{H}}$. Then
$\mathrm{ker}(\theta)=\bar{\mathbb{H}}$ by
Lemma~\ref{lemma:A4-non-trivial-action}.

The homomorphism $\theta$ induces an~isomorphism
$\tau\colon\mathbb{M}\to\SL_2(\mathbb{F}_{5})$.

Let $F\subset\mathbb{M}$ be a~subgroup such that $\beta(F)=\beta(G)\cong
2.\A_4$. Then $G=\mathbb{H}\rtimes F$.

Suppose that the~group $\bar{\mathbb{H}}$ contains a~proper
non-trivial subgroup that is a~normal subgroup of the~group
$\bar{G}$. Let us consider $\bar{\mathbb{H}}$ as a~two-dimensional
vector space over~\mbox{$\mathbb{F}_5$}. Then
$\mathbb{F}_{5}^{2}\cong\bar{\mathbb{H}}=V_0\oplus V_1$, where
$V_0$ and $V_1$ are one-dimensional $\tau(F)$-invariant subspaces,
since $|2.\A_4|=24$ is coprime~to~$5$.

By Lemma~\ref{lemma:A4-non-trivial-action},  the~homomorphism $\tau$ induces
a~monomorphism
$$
F\longrightarrow\GL_1\Big(\mathbb{F}_5\Big)\times\GL_1\Big(\mathbb{F}_5\Big)\cong\mathbb{Z}_4\times\mathbb{Z}_4,
$$
which implies that $F$ is an~abelian group, which is not the~case.
\end{proof}

\begin{lemma}
\label{lemma:A4-in-HM-normal-subgroups} The group $\bar{G}$ does
not contain proper normal subgroups not containing
$\bar{\mathbb{H}}$.
\end{lemma}

\begin{proof}
Suppose that $\bar{G}$ contains a~normal subgroup
$\bar{G}^{\prime}$ such that
$\bar{\mathbb{H}}\not\subseteq\bar{G}^{\prime}$. Then the
intersection $\bar{G}^{\prime}\cap\bar{\mathbb{H}}$ consists of
the identity element in $G$ by Lemma~\ref{lemma:A4-active-action}.
Hence
$$
\bar{G}^{\prime}\cong\beta\big(\bar{G}^{\prime}\big)\subseteq\beta\big(\bar{G}\big)\cong 2.\A_{4},%
$$
which implies that $\bar{G}^{\prime}$ is isomorphic to a~normal
subgroup of the~group $2.\A_4$.

Let $\bar{Z}$ be the~center of $\bar{G}^{\prime}$. Then $\bar{Z}$
is a~normal subgroup of the~group $\bar{G}$. Thus, we have
$\bar{Z}\cong\mathbb{Z}_{2}$ by
Lemma~\ref{lemma:normal-subgroups-in-A4}. Hence $\bar{Z}$ is
contained in the~center of $\bar{G}$, which contradicts
Lemma~\ref{lemma:A4-non-trivial-action}.
\end{proof}

\begin{lemma}
\label{lemma:A4-in-HM-elliptic} Let $E$ be a~smooth irreducible
curve of genus $g\le 8$.
Then there is no monomorphism $\bar{G}\to \mathrm{Aut}(E)$.
\end{lemma}

\begin{proof}
By classification of finite subgroups in
$\PGL_2(\mathbb{C})$ the case~\mbox{$g=0$} is impossible.
The cases~\mbox{$2\le g\le 8$} are
impossible by Theorem~\ref{theorem:Hurwitz-bound}.
Therefore, we may assume that $E$ is an elliptic curve.

Let us consider $E$ as an abelian group. Then there is an exact
sequence
$$
\xymatrix{&1\ar@{->}[rr]&&E\ar@{->}[rr]^{\iota}&&\mathrm{Aut}\big(E\big)\ar@{->}[rr]^{\upsilon}&&\mathbb{Z}_{n}\ar@{->}[rr]&&1}%
$$
for some $n\in\{2,4,6\}$.

Suppose that there is a~monomorphism $\theta\colon\bar{G}\to
\mathrm{Aut}(E)$. Then $\theta(\bar{\mathbb{H}})\subset\iota(E)$,
because $\iota(E)$ contains all the~elements of $\mathrm{Aut}(E)$
of order $5$.

Let $g$ be any element of $\bar{G}$ such that
$\theta(g)\in\iota(E)$. Then
$\theta(g)\theta(h)=\theta(h)\theta(g)$ for every
$h\in\bar{\mathbb{H}}$, because $\iota(E)$ is an abelian group,
and thus $g\in\bar{\mathbb{H}}$ by
Lemma~\ref{lemma:A4-non-trivial-action}. Hence
$\theta(\bar{G})\cap\iota(E)=\theta(\bar{\mathbb{H}})$, which
implies that $\upsilon(\bar{G})\cong\beta(\bar{G})\cong
2.\A_4$, which is absurd.
\end{proof}

The main purpose of this section is to prove the~following result.

\begin{theorem}
\label{theorem:HM-splitting} Let $E$ be~a~smooth elliptic curve.
Then there is no exact sequence of groups
\begin{equation}
\label{equation:HM-splittinng-sequence}
\xymatrix{&1\ar@{->}[rr]&&G^{\prime}\ar@{->}[rr]^{\iota}&&\bar{G}\ar@{->}[rr]^{\upsilon}&&G^{\prime\prime}\ar@{->}[rr]&&1,}%
\end{equation}
where $G^{\prime}$ and $G^{\prime\prime}$ are subgroups of
the~groups $\mathrm{Aut}(\mathbb{P}^1)$ and $\mathrm{Aut}(E)$,
respectively.
\end{theorem}

\begin{proof}
Suppose that the~exact sequence of
groups~\ref{equation:HM-splittinng-sequence} does exist. Then
$\iota$ is not an~isomorphism, because the~group
$\mathrm{Aut}(\mathbb{P}^1)$ does not contain subgroups isomorphic
to $\bar{G}$. The monomorphism $\upsilon$ is not an~isomorphism by
Lemma~\ref{lemma:A4-in-HM-elliptic}. Then
\mbox{$\bar{\mathbb{H}}\subset\iota(G^{\prime})$} by
Lemma~\ref{lemma:A4-in-HM-normal-subgroups}. But
$\mathrm{Aut}(\mathbb{P}^1)$
contains no subgroups isomorphic to $\bar{\mathbb{H}}$, which is
a~contradiction.
\end{proof}

\begin{corollary}
\label{corollary:elliptic-surface} There is no~monomorphism
$\bar{G}\to\mathrm{Bir}(E\times\mathbb{P}^{1})$, where $E$ is
a smooth elliptic curve.
\end{corollary}

We believe that there is a~simpler proof of Theorem~\ref{theorem:HM-splitting}.

\end{document}